\documentclass[11pt]{article}
\usepackage{amsmath}
\usepackage{amssymb}

{\catcode `\@=11 \global\let\AddToReset=\@addtoreset}
\AddToReset{equation}{section}

\usepackage{epsfig}
\usepackage{amsmath}
\usepackage{amssymb}
\usepackage{graphicx}
\usepackage{color}

\usepackage[latin1]{inputenc}

\newtheorem{theorem}{Theorem}[section]

\newtheorem{corollary}{\bf Corollary}[section]

\newtheorem{proposition}{Proposition}[section]
\newtheorem{@definition}{\sc Definition}[section]

\newtheorem{@remark}{\sc Remark}[section]
\newenvironment{remark}{\begin{@remark}\rm}{\end{@remark}}
\newtheorem{@example}{\sc Example}[section]

\newcommand{\beqn}{\begin{displaymath}}
\newcommand{\eeqn}{\end{displaymath}}
\newcommand{\beq}{\begin{equation}}  
\newcommand{\eeq}{\end{equation}}

\def\mathsf{\bf}
\def\N{\mathbb{N}}

\def\R{\mathbb{R}}
\def\Z{\mathbb{Z}}

\def\i{{\mathrm i}}
\def\d{{\mathrm d}}
\def\e{{\mathrm e}}

\def\E{\mathrm E}
\def\P{\mathrm P}

\def\text{\mbox}

\def\1{{\bf 1}}

\def\vep{\varepsilon}

\newcommand{\nn}{\nonumber}
\newcommand{\noi}{\noindent}

\def\eq2{
\stackrel{\small \rm mod \,2}{=}}

\def\n2{
\stackrel{\small \rm mod \,2}{\neq}}

\def\limd{\renewcommand{\arraystretch}{0.5}
\begin{array}[t]{c}
\stackrel{\rm d}{\longrightarrow} \\
\end{array}\renewcommand{\arraystretch}{1}}

\def\limfdd{\renewcommand{\arraystretch}{0.5}
\begin{array}[t]{c}
\stackrel{\rm fdd}{\longrightarrow} \\
\end{array}\renewcommand{\arraystretch}{1}}

\def\eqfdd{\renewcommand{\arraystretch}{0.5}
\begin{array}[t]{c}
\stackrel{\rm fdd}{=} \\
\end{array}\renewcommand{\arraystretch}{1}}

\def\neqfdd{\renewcommand{\arraystretch}{0.5}
\begin{array}[t]{c}
\stackrel{\rm fdd}{\neq} \\
\end{array}\renewcommand{\arraystretch}{1}}

\def\eqd{\renewcommand{\arraystretch}{0.5}
\begin{array}[t]{c}
\stackrel{\rm d}{=} \\
\end{array}\renewcommand{\arraystretch}{1}}

\newtheorem{thm}{Theorem}[section]

\newtheorem{lem}[thm]{Lemma}

\def\vep{\varepsilon}

\begin{document}


\title{Scaling  transition
for nonlinear random fields \\ with long-range dependence
}


\author{Vytaut\.e Pilipauskait\.e$^{1,2}$ and
Donatas Surgailis$^1$
\\
\small
$^1$ Vilnius University, 
Akademijos 4, LT-08663 Vilnius, Lithuania\\
\small
$^2$ Universit\'{e} de Nantes, 1, quai de Tourville BP, Nantes  44313, France
}
\maketitle

\begin{quote}

{\bf Abstract.}  We obtain a complete description of anisotropic scaling limits and
the existence of scaling transition
for nonlinear functions (Appell polynomials) of stationary linear random fields on $\Z^2$ with moving average coefficients
decaying at possibly different rate in the horizontal  and vertical direction.
The paper extends recent results on scaling transition  for
linear random fields in \cite{ps2014},  \cite{ps2015}.

\end{quote}

\smallskip
{\small

\noi {\it Keywords:} scaling  transition; anisotropic long-range dependence;
fractionally integrated random field;
Appell polynomials; multiple It\^o-Wiener integral;
fractional Brownian sheet

}

\vskip.7cm

\section{Introduction}

\cite{ps2014} introduced the notion of scaling transition for stationary random field (RF)
$X = \{ X(t,s); (t,s) \in \Z^2\}$ on  $\Z^2$ in terms of partial sums limits
\begin{equation} \label{Xsum01}
D^{-1}_{\lambda,\gamma} \sum_{(t,s) \in K_{[\lambda x, \lambda^{\gamma}y]}} X(t,s) \ \limfdd \ V_\gamma (x,y), \quad (x,y) \in \R^2_+, \quad \lambda \to \infty, \quad \gamma >0
\end{equation}
where $D_{\lambda,\gamma} \to \infty $ is normalization and
$K_{[\lambda x, \lambda^{\gamma}y]} := \{ (t,s)  \in \Z^2: 1\le t \le \lambda x, 1\le s \le \lambda^{\gamma}y \} $ is a family of rectangles
whose sides grow at possibly different rate  $O(\lambda)$ and $O(\lambda^{\gamma})$ and  $\gamma >0$ is {\it arbitrary}.
See the end of this section for all unexplained notation.
RF $X$ is said to exhibit {\it scaling transition} at $\gamma_0 >0$
if the limit RFs $V_\gamma \equiv V^X_\gamma $ in \eqref{Xsum01}
do not depend on $\gamma $ for  $\gamma >\gamma_0$ and $\gamma < \gamma_0$
and are different up to a multiplicative constant, viz.,
\begin{equation}
 V^X_\gamma \eqfdd V^X_+ \ \ (\forall \gamma > \gamma_0), \qquad
V^X_\gamma \eqfdd V^X_- \ \ (\forall \gamma < \gamma_0), \qquad V^X_+ \neqfdd a V^X_- \ \ (\forall a >0).
\end{equation}
In such case, RF  $V^X_{\gamma_0}$ 
is called the {\it well-balanced}
while RFs $V^X_+$ and  $V^X_-$ the {\it unbalanced} scaling limits of $X$.

It appears that scaling transition is a new and general feature of spatial dependence which  occurs for many
isotropic and anisotropic RF on $\Z^2$ with long-range dependence (LRD). It was established
for a class of aggregated $\alpha$-stable autoregressive models \cite{ps2014},
a class of Gaussian LRD RFs  \cite{ps2015}, and some RFs arising by aggregation of network traffic and random-coefficient time series
models in telecommunications and economics; see  \cite{gaig2003},  \cite{miko2002},  \cite{pils2014}, \cite{pils2015},
also  (\cite{ps2014}, Remark 2.3).
The unbalanced limits $V_\pm$ in the these studies have a very special dependence
structure  (either independent or invariant rectangular increments along one of the coordinate axes) and coincide
in the Gaussian case 
with a fractional Brownian sheet (FBS) $B_{H_1,H_2} $
with one of the two parameters   $H_1, H_2 \in (0,1]$ equal to $1/2$ or $1$.


The above mentioned works deal with linear RF models written as sums (stochastic integrals) w.r.t. i.i.d. `noise'. It is well-known
that nonlinear RFs can display quite complicated nongaussian scaling behavior. See  Dobrushin and Major \cite{dobmaj1979},
also \cite{AnhLRM2012}, \cite{avr1987}, \cite{gir1985}, \cite{gir2012},  \cite{hoh1997},
\cite{lav2007}, \cite{leo1999}, 
 \cite{sur1982}, \cite{taq1979} and the references therein. 

The present paper  establishes the existence of scaling transition for a class of
nonlinear  subordinated RFs:
\begin{equation}\label{Xnonlin}
X(t,s) =  G(Y(t,s)),
\end{equation}
where $Y = \{ Y(t,s), (t,s) \in \Z^2 \} $ is a stationary linear LRD RF in \eqref{Xlin} and $G(x)= A_k(x), x \in \R $ is the
Appell polynomial  of degree $k \ge 1$ (see Sec.~2 for the definition)
with $\E G^2(Y(0,0)) < \infty, \, \E G(Y(t,s)) = 0 $.
The (underlying) RF $Y$ is written as  a moving-average
\begin{equation}\label{Xlin}
Y(t,s) \ = \ \sum_{(u,v) \in \Z^2} a(t-u, s-v) \vep(u,v), \qquad (t,s) \in \Z^2,
\end{equation}
in a standardized i.i.d. sequence $\{ \vep(u,v); (u,v) \in \Z^2\}$ with deterministic  
moving-average coefficients such that
\begin{equation}\label{acoef}
a(t,s) \ \sim \ {\rm const} \, (|t|^2 +  |s|^{2q_2/q_1})^{-q_1/2},  \qquad |t|+|s| \to \infty,
\end{equation}
where parameters $q_1,q_2 >0$ satisfy
\begin{equation}\label{qLRD}
1 <  Q := \frac{1}{q_1} + \frac{1}{ q_2}  < 2.
\end{equation}
In Theorems \ref{thm1}-\ref{thm5} below, the moving-average coefficients $a(t,s)$ 
may take a more general form in \eqref{acoefL} 
including an `angular function'. 
Condition $Q < 2 $ guarantees that $\sum_{(t,s) \in \Z} a(t,s)^2 < \infty $ or $Y$ in \eqref{Xlin} is well-defined,
while $Q> 1$ implies that $\sum_{(t,s) \in \Z } |a(t,s)| = \infty $ (in other words, that RF $Y$ is LRD).
Note $a(t,0) = O(|t|^{-q_1}), \, a(0,s) = O(|s|^{-q_2})$ decay at a different rate when $q_1 \ne q_2$ in which case $Y$ 
exhibits strong anisotropy. The form of moving-average coefficients in \eqref{acoef} implies a similar behavior
of the covariance function
$r_Y(t,s) := \E Y(0,0) Y(t,s) =  \sum_{(u,v) \in \Z^2} a(u,v) a(t+u,s+v) $, namely
\begin{eqnarray}\label{covYdef}
r_Y(t,s)&\sim&{\rm const}\,
(|t|^{2} +|s|^{2p_2/p_1})^{-p_1/2}, \qquad |t|+|s| \to \infty,
\end{eqnarray}
where
\begin{equation}\label{pi}
p_i := q_i(2- Q), \qquad i=1,2.
\end{equation}
Note $p_1/p_2 = q_1/q_2$ and the 1-1 correspondence between $(q_1,q_2)$ and $(p_1,p_2)$:
\begin{equation}\label{qi}
q_i := \frac{p_i}{2}(1+ P), \qquad i=1,2, \quad \text{where} \quad  P := \frac{1}{p_1} + \frac{1}{p_2}.
\end{equation}
\eqref{covYdef} implies that for any integer
$k\ge 1 $ and $P \not\in \N$
\begin{eqnarray}\label{QP}
\sum_{(t,s) \in \Z^2} |r_Y(t,s)|^k = \infty  \qquad \Longleftrightarrow \qquad   1\le k < P,
\end{eqnarray}
see Proposition~\ref{prop1}. In the case when
$Y$ in \eqref{Xlin} is Gaussian RF, $r^k_{Y}(t,s)$ coincides with the covariance of the
$k$th Hermite polynomial $H_k(Y(t,s)) $  of $Y$
and the (nonlinear) subordinated RF  $X =  H_k (Y)$ is LRD if condition \eqref{QP} holds.
A similar result is true for nongaussian moving-average RF $Y$ in \eqref{Xlin} and
Hermite polynomial $H_k$ replaced by Appell polynomial $A_k$.

The following summary  describes the main results of this paper.

\begin{itemize}

\item[(R1)] Subordinated RFs $X = A_k(Y),  1\le k < P$
exhibit scaling transition at the same point $\gamma_0 := p_1/p_2 = q_1/q_2 $ independent of $k$.

\item[(R2)] The well-balanced scaling limit $V^X_{\gamma_0}$ of $X =  A_k(Y)$
is non-gaussian unless $k = 1$ and is given by a $k$-tuple  It\^o-Wiener integral.

\item[(R3)]
Unbalanced scaling limits $  V^X_+ = V^X_\gamma, \gamma >\gamma_0$ of $X =  A_k(Y)$
agree with FBS $B_{H^+_{1k},1/2}$ with Hurst parameter  $H^+_{1k} \in (1/2,1)  $ if $k p_2 > 1 $, and with a `generalized Hermite slide'
$V^X_{+}(x,y) = x Z^+_{k}(y)$ if  $k p_2 < 1 $, where $Z^+_{k}$ is a 
self-similar process
written as a $k$-tuple It\^o-Wiener integral.
A similar fact holds for unbalanced limits $  V^X_- = V^X_\gamma, \gamma < \gamma_0$.

\item[(R4)]\
For $k > P$, RF $X = A_k(Y)$  does not exhibit  scaling transition and all
scaling limits $V^X_\gamma, \gamma >0 $ agree with Brownian sheet $B_{1/2,1/2}$.

\item[(R5)]
In the case of Gaussian underlying RF $Y$ in \eqref{Xlin}, the above conclusions hold for general nonlinear  function $G$ in
\eqref{Xnonlin} and  $k$ equal to  the Hermite rank of $G$.

\end{itemize}

The above list contains several new noncentral and central limit results.
(R2), (R4) and (R5) are new in the `anisotropic' case $p_1 \ne p_2$ while
(R3) is new  even for linear RF $X = A_1(Y) = Y$ (see Remark~\ref{hermite} concerning the terminology in (R3)).
Similarly as in the case of linear models (see \cite{ps2014}, \cite{ps2015}),
unbalanced limits in (R3) have either independent or completely dependent increments along one of the coordinate axes.
According to (R3), the sample mean of  nonlinear LRD RF $X = A_k(Y), 1 < k < P$ on rectangles $K_{[\lambda, \lambda^\gamma]}, \gamma \ne \gamma_0$
may have  gaussian or nongaussian limit distribution depending on $k, \gamma$ and parameters $p_1, p_2$, moreover, in both cases
the variance of the sum $\sum_{(t,s)\in K_{[\lambda, \lambda^\gamma]} } X(t,s) $ grows faster than $\lambda^{1+ \gamma}$, or
the number of   summands.
The dichotomy of the limit distribution  in (R3)
is related to the presence or absence of the vertical/horizontal LRD property of $X$, see Remark~\ref{remLRD}.
We also note that
our proofs of the central limit results in (R3) and (R4) use rather simple approximation by $m$-dependent r.v.'s
and do not require a combinatorial argument or Malliavin's calculus as in \cite{breu1983}, \cite{nua2005}
and other papers.

The paper is organized as follows.  Sec.~2 provides the precise assumptions
on RFs $Y$ and $X$ and some known properties of Appell polynomials.
Sec.~3 contains formulations of the main results (Theorems~\ref{thm1}-\ref{thm5}) as described
in (R1)-(R5) above. Sec.~4 provides
two examples of linear fractionally integrated RFs satisfying the assumptions in Sec.~2.
Sec.~5 discusses some properties of generalized homogeneous functions and their convolutions
used to prove the results.
Sec.~6 discusses the asymptotic form of the covariance
function and  
the asymptotics of the variance of anisotropic
partial sums of subordinated RF $ X = A_k(Y)$. All proofs are collected in Sec.~7 and 8.

\smallskip

{\it Notation.} In this paper, $\limd, \limfdd, \eqd, \eqfdd$ denote the weak convergence and equality of (finite dimensional) distributions. $C$ stands for a generic positive constant which may assume different values at various locations and whose precise value has no importance.
$\R_+ := (0,\infty), \, \R^2_+ := (0,\infty)^2, \, \R^2_0 := \R^2 \setminus \{(0,0)\}, \, \Z_+ := \{0,1, \cdots \}, \, \N_+ := \{1,2, \cdots \}, \,
\Z^{\bullet 2k}  :=  \{ ((u_1,v_1), \cdots, (u_k,v_k)) \in \Z^{2k}: (u_i,v_i) \ne (u_j,v_j), 1 \le  i < j \le k \}, k \in \N_+, \,
|t|_+ :=  |t| \vee 1 \, (t \in \Z) $.

\section{Assumptions and preliminaries}

\noi {\bf Assumption (A1)} \ $\{\vep, \vep(t,s), (t,s) \in \Z^2 \}$ is an i.i.d. sequence with $\E \vep = 0, \E \vep^2 = 1 $

\medskip

\noi  {\bf Assumption (A2)} \ $ Y = \{ Y(t,s), (t,s) \in \Z^2 \} $ is a moving-average RF in \eqref{Xlin} with
coefficients
\begin{equation}\label{acoefL}
a(t,s)\  = \  \frac{1}{ (|t|^2 +  |s|^{2q_2/q_1})^{q_1/2}}\Big(L_0 \big(\frac{t}{(|t|^{2} +|s|^{2q_2/q_1})^{1/2}}\big) + o(1)\Big),
\qquad |t|+|s| \to  \infty,
\end{equation}
where $q_i>0, i=1,2 $ satisfy $Q = \sum_{i=1}^2 q_i^{-1} \in (1,2) $ (see \eqref{qLRD})
and $L_0(u) \ge 0, u \in [-1,1]$ is a bounded piece-wise continuous function on $[-1,1]$.

\medskip

We refer to $L_0$ in \eqref{acoefL} as the {\it angular function}.
Particularly, for $q_1 = q_2 $,  $\rho = (|t|^2 +  |s|^{2})^{1/2} $ and
$\arccos (t/\rho) $ are the polar coordinates of $(t,s) \in \R^2$.
Assumptions (A1)-(A2) imply $\E Y(0,0)^2 =  \sum_{(t,s) \in \Z^2} a(t,s)^2  $  $ < \infty $ and hence RF $Y$
in \eqref{Xlin} is well-defined and stationary, with zero mean $\E Y(t,s) = 0$. Moreover,
if $\E |\vep|^\alpha < \infty $ for some $\alpha > 2 $ then $\E |Y(t,s)|^\alpha < \infty $ follows
by  Rosenthal's inequality; see e.g. (\cite{gir2012}, Corollary~2.5.1).

Given a r.v. $\xi $ with $\E |\xi|^k < \infty, k \in \Z_+$, the $k$th Appell polynomial
$A_k(x)$ relative to the distribution of $\xi$ is defined by  $A_k(x) := (-\i)^k \d^k (\e^{\i u x}/\E \e^{\i u \xi})/\d u^k\big|_{u=0}$.
See \cite{avr1987}, \cite{gir2012} for
various properties of Appell polynomials.
In as follows, $A_k (\xi) $ stands for the r.v. obtained by substituting
$x= \xi$  in the Appell polynomial $A_k(x)$ relative to the distribution of $\xi $. Particularly, if $\E \xi = 0$ then
$A_1(\xi) = \xi,  \,
A_2(\xi) = \xi^2 - \E \xi^2, \, A_3(\xi) = \xi^3 - 3 \xi  \E \xi^2 - \E \xi^3 $ etc.
For standard normal $\xi \sim N(0,1)$ the Appell polynomials
$A_k(\xi) = H_k(\xi) = (-\i)^k \d^k\e^{\i u \xi  + u^2/2}/\d u^k \big|_{u=0} $ agree with the Hermite polynomials.

\medskip

\noi  {\bf Assumption {\bf (A3)}$_k$} \ For $k \in \N_+ $,  $\E |\vep|^{2k} < \infty $ and
\begin{equation} \label{AY}
X = \{ X(t,s) := A_k(Y(t,s)), (t,s) \in \Z^2 \},
\end{equation}
where $A_k $ is the $k$th Appell polynomial relative to the (marginal) distribution of $Y(t,s)$ in
\eqref{Xlin}.

\medskip

We also use the representation
of \eqref{AY} via {\it Wick products} of noise variables (see \cite{gir2012}, Ch.~14):
\begin{eqnarray}\label{Awick}
A_k(Y(t,s)) \ = \
\sum_{(u,v)_k \in \Z^{2k}}
a(t-u_1,s-v_1) \cdots a(t-u_k,s-v_k) :\!\vep(u_1,v_1) \cdots \vep(u_k,v_k)\!:.
\end{eqnarray}
By definition, for {\it mutually distinct} points
$(u_j, v_j) \ne (u_{j'},v_{j'}) \ (j \ne j', 1 \le j,j' \le i)$  the Wick product $:\!\vep(u_1,v_1)^{k_1} \cdots \vep(u_i,v_i)^{k_i}\!: \, = \prod_{j=1}^i A_{k_j}(\vep(u_j,v_j)) $ equals the product
of {\it independent} r.v.'s $ A_{k_j}(\vep(u_j,v_j)), 1 \le j \le i$. \eqref{Awick} leads to the decomposition
of \eqref{AY} into the `off-diagonal' and `diagonal' parts:
\begin{equation}\label{Akdec}
A_k(Y(t,s)) = Y^{\bullet k}(t,s) + {\cal Z}(t,s),
\end{equation}
where
\begin{eqnarray}\label{Ykstar}
Y^{\bullet k}(t,s)&:=&{\sum_{(u,v)_k}}^{\!\!\!\bf \bullet}  a(t-u_1,s-v_1) \cdots a(t-u_k,s-v_k) \vep(u_1,v_1) \cdots \vep(u_k,v_k)
\end{eqnarray}
and the sum $\sum_{(u,v)_k}^\bullet$ is taken over all $(u,v)_k = ((u_1,v_1), \cdots, (u_k,v_k)) \in \Z^{2k}$ such that
$(u_i,v_i) \ne (u_j,v_j) \ (i \ne j)$ (the set of such $(u,v)_k \in \Z^{2k} $ will be denoted by $\Z^{\bullet 2k} $).
By definition, the `diagonal' part  ${\cal Z}(t,s)$ in  \eqref{Akdec} is given by the r.h.s. of \eqref{Awick} with
$(u,v)_k \in \Z^{2k}$ replaced by $(u,v)_k \in \Z^{2k}\setminus \Z^{\bullet 2k}. $ In most of our limit results,
${\cal Z}(t,s)$ is negligible and  $Y^{\bullet k}(t,s)$ is the main term which is easier to handle compared to $A_k(Y(t,s))$
in \eqref{Akdec}. We also note that limit distributions of partial sums of `off-diagonal' polynomial forms in i.i.d. r.v.'s
were studied in \cite{sur1982}, \cite{gir2012}, \cite{bai2014} and other works.

\medskip

\noi  {\bf Assumption {\bf (A4)}$_k$} \ $ \vep(0,0) \eqd Z $ and $ Y(0,0) \eqd Z$ have  standard normal distribution
$Z \sim N(0,1)$
and $ X(t,s) = G(Y(t,s))$, where $G = G(x), x \in \R$ is a measurable function with $\E G(Z)^2 < \infty, \E G(Z) = 0$ and
Hermite rank $k \ge 1 $.

\medskip

Assumptions (A1), (A2) and (A4)$_k$ imply that $Y$ in \eqref{Xlin} is a Gaussian RF.
As noted above, under Assumption (A4)$_k$ Appell polynomials $A_k(x)$ coincide with Hermite polynomials $H_k(x)$. Recall that
the Hermite rank of a measurable function $G: \R \to \R$ with $\E G(Z)^2 < \infty$ is defined as the index $k$ of the lowest nonzero coefficient $c_j$
in the Hermite expansion of $G$, viz., $G(x) = \sum_{j=k}^\infty c_j H_j(x)/j! $ where $c_k \ne 0$.

Let  $L^2({\R}^{2k})$ denote the Hilbert space of real-valued functions $h = h((u,v)_k), (u,v)_k = (u_1,v_1,
\cdots, u_k,v_k) \in \R^{2k}$ with finite norm $\|h\|_k :=  \{ \int_{\R^{2k}} h^2((u,v)_k)\d (u,v)_k\}^{1/2}, $  $ \d (u,v)_k =  \d u_1 \d v_1 \cdots
\d u_k \d v_k. $  Let  $W = \{ W(\d u, \d v), $  $ (u,v) \in \R^2 \}$ denote a real-valued Gaussian white noise
with zero mean and variance $\E W(\d u, \d v)^2 = \d u \d v$.
For any $h \in L^2({\R}^{2k})$
the $k$-tuple It\^o-Wiener integral
$\int_{\R^{2k}} h((u,v)_k)\d^k W =  \int_{\R^{2k}} h(u_1,v_1, \cdots, u_k,v_k) W(\d u_1,\d v_1) \cdots $  $ W(\d u_k,\d v_k)$
is well-defined
and satisfies  $\E \int_{\R^{2k}} h((u,v)_k)\d^k W = 0, \,
\E \big(\int_{\R^{2k}} h((u,v)_k)\d^k W\big)^2 \le k! \|h\|^2_k; $
see  e.g. \cite{gir2012}.

\section{Main results}

Recall the definitions $p_i, P$ in \eqref{pi}, \eqref{qi}; $\gamma_0 = q_1/q_2 = p_1/p_2$. Denote
\begin{equation}
S^X_{\lambda,\gamma}(x,y)  :=  \sum_{(t,s) \in K_{[\lambda x,\lambda^{\gamma} y]}}  X(t,s), \qquad S^X_{\lambda,\gamma} :=  S^X_{\lambda,\gamma}(1,1).
\end{equation}
Consider a RF
\begin{equation} \label{V00}
V_{k,\gamma_0}(x,y)\ := \  \int_{\R^{2k}}h(x,y;(u,v)_k)\d^k W, \quad  (x,y) \in \R^2_+,
\end{equation}
where (c.f. \eqref{acoefL})
\begin{eqnarray} \label{V0}
h(x,y;(u,v)_k)&:=&\int_{(0,x]\times (0,y]} \prod_{\ell =1}^k a_\infty(t-u_\ell, s-v_\ell) \d t \d s, \quad \text{where} \\
a_\infty (t,s)&:=&(|t|^2 +  |s|^{2q_2/q_1})^{-q_1/2} L_0 \big(t/(|t|^{2} +|s|^{2q_2/q_1})^{1/2}\big), \quad (t,s) \in \R^2. \nn
\end{eqnarray}

\begin{theorem} \label{thm1} (i) The RF $V_{k,\gamma_0}$ in \eqref{V00}
is well-defined for $1\le k < P $ as It\^o-Wiener
stochastic integral  and has zero mean
$\E V_{k,\gamma_0}(x,y) = 0$  and finite variance
$\E V^2_{k,\gamma_0}(x,y) = k! \|h(x,y; \cdot)\|^2_k$.
Moreover, RF $V_{k,\gamma_0}$ has stationary rectangular increments and satisfies
the OSRF property:
\begin{equation} \label{Vss}
\{ V_{k,\gamma_0}(\lambda x, \lambda^{\gamma_0} y), (x,y) \in \R^2_+   \}
\eqfdd \{\lambda^{H(\gamma_0)} V_{k,\gamma_0}(x,y), (x,y) \in \R^2_+   \}, \qquad \forall  \lambda >0,
\end{equation}
where $H(\gamma_0) := 1 + \gamma_0 - kp_1/2.$

\noi (ii) Let RFs $Y$ and $ X = A_k(Y)$ satisfy Assumptions (A1), (A2) and (A3)$_k$, $1\le k < P$.
Then
\begin{equation}\label{varY}
{\rm Var}(S^X_{\lambda,\gamma_0})
\  \sim\  c(\gamma_0)\lambda^{2H(\gamma_0)}, \qquad c(\gamma_0) := \|h(1,1;\cdot)\|^2_k
\end{equation}
and
\begin{equation}\label{limY}
n^{-H(\gamma_0)} S^X_{\lambda,\gamma_0}(x,y)
\ \limfdd\  V_{\gamma_0}(x,y).
\end{equation}
\end{theorem}

Next, we discuss the case $  k < P, \ \gamma \neq \gamma_0.$  This case is split into four subcases: (c1):  $\gamma > \gamma_0,  k > 1/p_2, $  (c2):  $\gamma > \gamma_0,  k < 1/p_2, $
(c3): $\gamma < \gamma_0,  k > 1/p_1, $  and (c4):  $\gamma < \gamma_0,  k < 1/p_1 $
(the `boundary'  cases $k = 1/p_i, i=1,2 $ are more delicate and omitted, see Remark \ref{kpi} below).
Cases (c3) and (c4) are symmetric to (c1) and (c2) and essentially follow by exchanging the coordinates $t$ and $s$.
Introduce random processes $Z^+_{k}$ and $Z^-_{k}$ with one-dimensional time:
\begin{eqnarray}\label{Zpm}
Z^+_{k}(y)&:=&\int_{\R^{2k}} h_+(y; (u,v)_k)  \d^k W,
\qquad
Z^-_{k}(x)\ :=\ \int_{\R^{2k}} h_-(x; (u,v)_k)  \d^k W, \qquad x, y \ge 0,
\end{eqnarray}
where
\begin{equation}\label{hpm}
h_+(y; (u,v)_k) \ := \ \int_0^y \prod_{i=1}^k  a_\infty (u_i, s-v_i) \d s, \qquad
h_-(y; (u,v)_k) \ := \ \int_0^x \prod_{i=1}^k  a_\infty (t-u_i, v_i) \d t,
\end{equation}
and $a_\infty(t,s) $ is defined in \eqref{V0}.

\begin{theorem} \label{thm2} (i) Processes  $Z^+_{k}$ and $Z^-_{k} $ in \eqref{Zpm} are well-defined
for $1\le k < 1/p_2$ and $1 \le k < 1/p_1$, respectively, as It\^o-Wiener stochastic integrals. They have zero mean, finite variance,
stationary increments and are self-similar with respective indices $H^+_{2k} := 1- kp_2/2 \in (1/2,1)$ and   $H^-_{1k} := 1- kp_1/2 \in (1/2,1)  $.

\smallskip

\noi (ii) Let RFs $Y$ and $ X = A_k(Y)$ satisfy Assumptions (A1), (A2) and (A3)$_k$, $1\le k < 1/p_2$.
Then for any $\gamma > \gamma_0$
\begin{equation}\label{varY2}
{\rm Var}(S^X_{\lambda,\gamma})
\ \sim\   c(\gamma) \lambda^{2H(\gamma)},
\end{equation}
where $H(\gamma) := 1 + \gamma H^+_{2k}$  and $c(\gamma) := \|h_+(1;\cdot)\|^2_k $.
Moreover,
\begin{equation}\label{limY2}
\lambda^{-H(\gamma)} S^X_{\lambda,\gamma}(x,y)
\ \limfdd\  x Z^+_{k}(y).
\end{equation}

\smallskip

\noi (iii) Let RFs $Y$ and $ X = A_k(Y)$ satisfy Assumptions (A1), (A2) and (A3)$_k$, $1\le k < 1/p_1$.
Then for any $\gamma < \gamma_0$
\begin{equation}\label{varY3}
{\rm Var}(S^X_{\lambda,\gamma})
\ \sim\  c(\gamma) \lambda^{2H(\gamma)}, 
\end{equation}
where $H(\gamma) := \gamma + H^-_{1k}$ and $c(\gamma) := \|h_-(1;\cdot)\|^2_k >0$.
Moreover,
\begin{equation}\label{limY3}
\lambda^{-H(\gamma)} S^X_{\lambda,\gamma}(x,y)
\ \limfdd\  y Z^-_{k}(x).
\end{equation}
\end{theorem}

\begin{remark} \label{hermite} Processes $Z^\pm_{k}$  in \eqref{Zpm} have a similar structure and properties to
generalized Hermite processes discussed in  \cite{bai2014} except that \eqref{Zpm} are defined
as $k$-tuple It\^o-Wiener integrals with respect to white noise in $\R^2 $ and not in $\R$ as in  \cite{bai2014}.
Following the terminology in \cite{pils2016}, RFs   $x Z^+_{k}(y)$ and  $y Z^-_{k}(x) $ may be called a
{\it generalized Hermite slide} since they represent a random surface `sliding linearly to $0$' along one of the coordinate on the plane
from a generalized Hermite process indexed by the other coordinate. In the Gaussian case $k=1$,
a  generalized Hermite slide agrees with a FBS $B_{H_1,H_2}$
where one of the two parameters $H_1, H_2$ equals 1.
Recall that a fractional Brownian sheet (FBS) $B_{H_1,H_2} =  \{B_{H_1,H_2}(x,y), (x,y) \in  \R^2_+\}$
with parameters  $0 < H_1, H_2 \le 1$ is a
Gaussian process 
with zero mean and covariance function
\begin{equation}\label{FBScov}
\E B_{H_1,H_2}(x_1,y_1) B_{H_1,H_2}(x_2,y_2) = (1/4) (x_1^{2H_1} + x_2^{2H_1} - |x_1-x_2|^{2H_1})
 (y_1^{2H_2} + y_2^{2H_2} - |y_1-y_2|^{2H_2}).
\end{equation}
\end{remark}

\begin{theorem} \label{thm3} (i) Let RFs $Y$ and $ X = A_k(Y)$ satisfy Assumptions (A1), (A2) and (A3)$_k$, $1/p_2 < k < P$.
Then for any $\gamma >\gamma_0$
\begin{equation}\label{varY4}
{\rm Var}(S^X_{\lambda,\gamma})
\ \sim\  c(\gamma) \lambda^{2H(\gamma)}, 
\end{equation}
where $H(\gamma) := H^+_{1k} + \gamma/2, \, H^+_{1k}:= 1  + \gamma_0/2 - kp_1/2 \in (1/2,1)$ and
$ c(\gamma)  := \int_{(0,1]^2 \times \R}
\big((a_{\infty} \star a_{\infty})(t_1-t_2,s)\big)^k  \d t_1 \d t_2  \d s  >0$.
Moreover,
\begin{equation}\label{limY4}
\lambda^{-H(\gamma)} S^X_{\lambda,\gamma}(x,y)
\ \limfdd\  c(\gamma)^{1/2} B_{H^+_{1k}, 1/2}(x,y).
\end{equation}

\smallskip

\noi (ii) Let RFs $Y$ and $ X = A_k(Y)$ satisfy Assumptions (A1), (A2) and (A3)$_k$, $1/p_1 < k < P$.
Then for any $\gamma <\gamma_0$
\begin{equation}\label{varY5}
{\rm Var}(S^X_{\lambda,\gamma})
\ \sim\  c(\gamma) \lambda^{2H(\gamma)}, 
\end{equation}
where $H(\gamma) := \gamma H^-_{2k} + 1/2,  \, H^-_{2k}:= 1  + 1/(2\gamma_0) - kp_2/2 \in (1/2,1)$
and $ c(\gamma)  := \int_{\R\times (0,1]^2 }
\big((a_{\infty} \star a_{\infty})(t,s_1-s_2)\big)^k  \d t \d s_1  \d s_2 >0 $.
Moreover,
\begin{equation}\label{limY5}
\lambda^{-H(\gamma)} S^X_{\lambda,\gamma}(x,y)
\ \limfdd\  c(\gamma)^{1/2} B_{1/2, H^-_{2k}}(x,y).
\end{equation}

\end{theorem}

\begin{remark} \label{kpi} Note $H^+_{1k} = 1 \ (kp_2=1)$ and $H^-_{2k} = 1 \ (kp_1=1)$. We expect that
the convergences \eqref{limY4} and \eqref{limY5} remain true (modulus a logarithmic correction of normalization)
in the `boundary' cases $ kp_2=1 $  and $kp_1=1 $ of Theorem \ref{thm3} (i) and (ii)
and the limit RFs in these cases agree with FBS $B_{1,1/2}$ or $B_{1/2,1}$, respectively, having both parameters
equal to 1 or 1/2.

\end{remark}

\smallskip

\noi The next theorem discusses the case $k > P$.

\begin{theorem} \label{thm4} Let RFs $Y$ and $ X = A_k(Y)$ satisfy Assumptions (A1), (A2) and (A3)$_k$ and $ k >P$.
Then for any $\gamma >0$
\begin{equation}\label{varY0}
{\rm Var}(S^X_{\lambda,\gamma})
\ \sim\  \sigma^2_X \lambda^{1+ \gamma},
\end{equation}
where $\sigma^2_X := \sum_{(t,s) \in \Z^2} {\rm Cov}(X(0,0), X(t,s)) \in (0, \infty) $.  Moreover,
\begin{eqnarray}\label{limY0}
\lambda^{-(1+\gamma)/2}S^X_{\lambda,\gamma}(x,y)
&\limfdd&\sigma_X B_{1/2,1/2}(x,y).
\end{eqnarray}
\end{theorem}

\smallskip

\noi Our last theorem extends the above results to general function $G$ having Hermite rank $k$ and
Gaussian underlying RF $Y$.

\begin{theorem} \label{thm5} Let $ X = G(Y)$ satisfy Assumption (A4)$_k$.
Assume w.l.g. that $G$ has  Hermite expansion  $G(x) = H_k(x) + \sum_{j=k+1}^\infty c_j H_j(x)/j! $.

\smallskip

\noi (i) Let $1 \le k < P$. Then RF $X$ satisfies all statements of Theorems~\ref{thm1}-\ref{thm3}.

\smallskip

\noi (ii) Let $k > P$. Then RF $X$ satisfies the statements of Theorem~\ref{thm4}.

\end{theorem}

According to Theorems~\ref{thm2}-\ref{thm3}, the unbalanced scaling limits $V^X_\pm $ of
RF  $X = A_k(Y)$ satisfying Assumptions (A1)-(A3)$_k$ are given by
\begin{equation}\label{Vpm}
V^X_+(x,y) \ = \  \begin{cases}
x Z^+_k(y),  &kp_2 < 1, \\
c_+^{1/2} B_{H^+_{1k},1/2}(x,y), &kp_2 > 1,
\end{cases}
\qquad
V^X_-(x,y) \ = \  \begin{cases}
y Z^-_k(x),  &kp_1 < 1, \\
c_-^{1/2}B_{1/2, H^-_{2k}}(x,y), &kp_1 > 1,
\end{cases}
\end{equation}
where $c_\pm \equiv c(\gamma) >0$ are given constants. The covariance functions of RFs $V^X_\pm$  in \eqref{Vpm}
agree (modulus a constant) with
the covariance of FBS $B_{H_1,H_2}$ where at least one of the two parameters $H_1, H_2$ equals $1$ or $1/2$, namely
$(H_1,H_2) = (1,H^+_{2k}) $ if $kp_2 < 1$, $=(H^+_{1k}, 1/2) $ if $ kp_2 > 1$ in the case of $V^X_+$, and
$(H_1,H_2) = (H^-_{1k},1) $ if $kp_1 < 1$, $=(1/2, H^-_{2k}) $ if $ kp_1 > 1$ in the case of $V^X_-$.
These facts and the explicit form of the covariance of FBS, see  \eqref{FBScov},
imply that $V_+ \neqfdd  a V_-  \, (\forall a >0)$, for any $k, p_1, p_2$ in Theorems~\ref{thm2}-\ref{thm3},
yielding the following corollary.

\begin{corollary} Let RF $ X = A_k(Y)$ satisfy Assumptions (A1), (A2) and (A3)$_k$, $1 \le k < P, kp_i \ne 1, $
$i=1,2.$  Then $X$ exhibits scaling transition at $\gamma_0 =  p_1/p_2$.
\end{corollary}

\section{Examples: fractionally integrated RFs}

In this section we present two examples of linear fractionally integrated random fields $Y$ in $\Z^2$
satisfying Assumptions (A1) and (A2).

\smallskip

\paragraph{Example 1. Isotropic fractionally integrated random field.} Introduce the (discrete) Laplace operator
$\Delta Y(t,s) := (1/4) \sum_{|u|+|v|=1} (Y(t+u,s+v) - Y(t,s))$
and a lattice isotropic fractionally integrated random field satisfying the equation:
\begin{equation}\label{Xfrac1}
(- \Delta)^d Y(t,s) \ = \  \vep(t,s),
\end{equation}
where $\{\vep(t,s), (t,s) \in \Z^2 \}$ are standard i.i.d. r.v.'s, $0< d < 1/2$ is the order of fractional integration,
$(1-z)^{d} =  \sum_{j=0}^\infty \psi_j(d) z^j,
\psi_j(d) := \Gamma(j-d)/\Gamma(j+1) \Gamma(- d). $ More explicitly,
\begin{eqnarray}\label{Xfrac11}
(- \Delta)^{d} Y(t,s)
&=&\sum_{j=0}^\infty \psi_j(d) (1+\Delta)^j Y(t,s)  \
=\  \sum_{(u,v) \in \Z^2} b(u,v) Y(t-u,s-v),
\end{eqnarray}
where $b(u,v):=\sum_{j=0}^\infty \psi_j(d) p_j(u,v) $
and $p_j (u,v) $ 
are $j$-step transition probabilities  of a symmetric  nearest-neighbor random walk $\{W_k, k=0,1,\cdots \}$
on $\Z^2$ with equal 1-step probabilities $P(W_1 = (u,v)|W_0 = (0,0)) =
1/4, |u| + |v| = 1$.
Note $\sum_{(u,v) \in \Z^2} |b(u,v)| = \sum_{j=0}^\infty |\psi_j(d)| < \infty, \, d > 0$ and therefore the l.h.s. of
(\ref{Xfrac11}) is well-defined
for any stationary random field $\{Y(t,s)\}$ with $\E |Y(0,0)| < \infty $.  As shown in
\cite{kou2016}, for $0< d < 1/2 $ a stationary solution of
(\ref{Xfrac11}) with zero-mean and finite  variance can be defined as a moving-average random field:
\begin{equation}\label{Xfrac13}
Y(t,s) = (- \Delta)^{-d} \vep(t,s) \ = \  \sum_{(u,v)\in \Z^2} a(u,v) \vep(t-u,s-v),
\end{equation}
with coefficients
\begin{equation}\label{bcoef1}
a(u,v) = \sum_{j=0}^\infty \psi_j(-d) p_j(u,v)
\end{equation}
satisfying $\sum_{(u,v) \in \Z^2} a(u,v)^2 < \infty. $  Moreover, RF $Y$ in \eqref{Xfrac13} has an explicit
spectral density $f(x,y) = (2\pi)^{-2} 2^{-2d}|(1-\cos x) +
(1-\cos y)|^{-2 d},  (x,y) \in [-\pi,\pi]^2$ which behaves as ${\rm const} \, (x^2 + y^2)^{-2d} $ as
$x^2 +  y^2 \to 0$.  According to (\cite{kou2016}, Proposition~5.1), the moving-average coefficients
in \eqref{bcoef1} satisfy the isotropic asymptotics:
\begin{equation*}
a(t,s) = (A + o(1))(t^2 + s^2)^{-(1- d)}, \quad t^2 + s^2 \to \infty,
\end{equation*}
where $A := \pi^{-1} \Gamma (1- d)/\Gamma (d)$ and hence Assumption (A2) with
$q_1 = q_2 = 2(1-d)\in (1,2),  Q = 1/(1-d) \in (1,2)$ and a constant angular function $L_0(z) = A, z \in [-1,1]$.

\paragraph{Example 2. Anisotropic fractionally integrated random field.} Consider the `discrete heat operator'
$\Delta_{1,2} Y(t,s) = Y(t,s) -  \theta Y(t-1, s) - \frac{1-\theta}{2} (Y(t-1,s+1) + Y(t-1,s-1)), 0< \theta < 1$
and a fractionally integrated random field satisfying
\begin{equation}\label{Xfrac2}
\Delta_{1,2}^d Y(t,s) \ = \  \vep(t,s),
\end{equation}
where $\{\vep(t,s)\}$ are as in \eqref{Xfrac1}. Similarly to \eqref{Xfrac13},
a stationary solution of  \eqref{Xfrac2} can be written
as a moving-average random field:
\begin{equation}\label{Xfrac21}
Y(t,s) = \Delta_{1,2}^{-d} \vep(t,s) \ = \  \sum_{(u,v)\in \Z_+ \times \Z} a(u,v) \vep(t-u,s-v),
\end{equation}
with coefficients
\begin{equation}\label{auv}
a(u,v) = \psi_u(-d) q_u(v)
\end{equation}
where  $q_u(v) $
are $u$-step transition probabilities  of a random walk $\{W_u, u=0,1,\cdots \}$
on $\Z$ with 1-step probabilities $P(W_1 = v|W_0 = 0) =
\theta $ if $v = 0$, $=  (1- \theta)/2$ if $v= \pm 1$. As shown in \cite{lls2014}, $\sum_{(u,v) \in \Z^2} a(t,s)^2 < \infty $
and the RF in \eqref{Xfrac21} is well-defined
for any $0< d < 3/4, \theta \in [0,1) $; moreover, the spectral density $f(x,y) $ of  \eqref{Xfrac21} is singular at the origin:
$f(x,y) \sim {\rm const}\, (x^2 + (1-\theta)^2 y^4/4)^{-d}, \,
(x,y) \to (0,0)$.

\begin{proposition}\label{propex2}  For any $ 0< d < 3/4,  0< \theta < 1 $
the coefficients in \eqref{auv} satisfy Assumption (A2) with
$q_1 = 3/2 -  d, q_2 = 2q_1 $ and a continuous angular function $L_0(z), z \in [-1,1]$ given by
\begin{equation}\label{L0}
L_0(z) \ = \ \begin{cases}
\frac{z^{d- 3/2}}{\Gamma (d) \sqrt{2 \pi (1-\theta)}}  \exp \big\{ - \frac{\sqrt{(1/z)^2 -1}}{2(1-\theta)}\big\},
&0 < z \le 1, \\
0, & -1 \le z \le 0.
\end{cases}
\end{equation}
\end{proposition}

\begin{remark} \label{negative}
\cite{bois2005}, \cite{guo2009} discussed fractionally integrated RFs  satisfying the equation
\begin{equation}\label{Xfrac3}
\Delta_1^{d_1} \Delta_2^{d_2}  Y(t,s) \ = \  \vep(t,s),
\end{equation}
where $\Delta_1 Y(t,s) :=  Y(t,s) - Y(t-1,s), \Delta_2 Y(t,s) := Y(t,s) - Y(t,s-1) $ are difference operators
and $0<d_1, d_2 < 1/2 $ are parameters. Stationary solution of \eqref{Xfrac3} is a moving-average RF in
$Y(t,s) = \sum_{(u,v) \in \Z_+^2} a(u,v)\vep(t-u,v-s)$ with coefficients
$a(u,v) := \psi_u(-d_1) \psi_v (-d_2).$  Following the proof of Theorem~\ref{thm1} one can show that for 
{\it any} $\gamma >0$  the (normalized) partial sums process of RF $Y$ in  \eqref{Xfrac3} tends to a FBS depending on $d_1, d_2$ only, viz.,
$\lambda^{-H_1 - \gamma H_2} S^Y_{\lambda,\gamma}(x,y) \limfdd c(d_1) c(d_2) B_{H_1,H_2}(x,y)  $, where
$H_i = d_i + 1/2 $ and $c(d_i) >0$ are some constants. See (\cite{ps2015}, Proposition~3.2) for related result. We conclude
that the fractionally integrated RF in \eqref{Xfrac3} featuring a `separation of LRD  along coordinate axes'
does not exhibit scaling transition in contrast to models in  \eqref{Xfrac1} and \eqref{Xfrac2}.

\end{remark}

\section{Properties of convolutions of generalized homogeneous functions}

For a given $\varpi >0$ denote
\begin{eqnarray}\label{rhodef}
&
\rho(t,s) :=  (|t|^2 + |s|^{2/\varpi})^{1/2}, \qquad
\rho_+(t,s) := 1 \vee \rho(t,s), \quad (t,s) \in \R^2.
\end{eqnarray}
Let $f(t,s) = \rho(t,s)^{-h} L(t/\rho(t,s))$, where $h \in \R $ and $L = L(z), z \in [-1,1]$ is an arbitrary measurable function, then
$f(t,s)$ satisfies the scaling property:
$f(\lambda t, \lambda^\varpi s) = \lambda^{-h} \rho(t,s), (t,s) \in \R^2_0 $
for each $\lambda >0$. Such functions are called generalized homogeneous functions (see \cite{han1972}).\\
We use the notation $[\phi_1 \star \phi_2](t,s) = \sum_{(u,v) \in \Z^2} \phi_1(u,v) \phi_2(t+u,s+v)$ for `discrete' convolution of
sequences $\phi_i = \{\phi_i(u,v), (u,v) \in \Z^2 \} \in L^2(\Z^2), i=1,2 $ and
$ (\psi_1 \star \psi_2)(t,s) = \int_{\R^2} \psi_1(u,v) \psi_2(t+u,s+v) \d u \d v $ for `usual'  convolution
of functions $\psi_i = \{\psi_i(u,v), (u,v) \in \R^2\}, i=1,2 $. Note the symmetry $\phi_i (t,s) = \phi_i(t,-s), i=1,2 $
implies the symmetry $[\phi_1 \star \phi_2](t,s) = [\phi_1 \star \phi_2](t,-s),   (\phi_1 \star \phi_2)(t,s)
= (\phi_1 \star \phi_2)(t,-s)$ of convolutions.
\\
Let $B_\delta(t,s) := \{ (u,v)\in \R^2: |t-u|+|s-v|\le \delta \}, \, B^c_\delta(t,s) := \R^2\setminus B_\delta(t,s)$.

\begin{proposition} \label{prop1} (i) For any $\delta >0, h >0$,
\begin{eqnarray}\label{rrho1}
&\int_{B_\delta(0,0)} \rho(t,s)^{-h} \d t \d s < \infty \ \  \Longleftrightarrow \ \  h < 1 + \varpi
\end{eqnarray}
and
\begin{eqnarray}\label{rrho2}
&\left.
\begin{array}{l}
\int_{B^c_\delta(0,0)} \rho(t,s)^{-h} \d t \d s < \infty \\
\sum_{(t,s)\in \Z^2} \rho_+(t,s)^{-h} < \infty
\end{array}\right\}
\ \ \Longleftrightarrow \  \ h > 1 + \varpi.
\end{eqnarray}

\smallskip

\noi (ii) Let $h_i >0, i=1,2, \, h_1 + h_2 > 1+ \varpi. $ Then there exists
$C>0$ such that for any $(t,s) \in \R^2_0 $
\begin{eqnarray}
(\rho^{-h_1} \star \rho^{-h_2})(t,s)&\le&C\rho(t,s)^{1 + \varpi - h_1 - h_2}, \quad \, h_i < 1 + \varpi, \, i=1,2, \label{prop11}  \\
(\rho_+^{-h_1} \star \rho^{-h_2})(t,s)&\le&C\rho_+(t,s)^{- h_2},  \hskip1.5cm h_2 < 1 + \varpi,  \, h_1 > 1 + \varpi,   \label{prop12}  \\
(\rho_+^{-h_1} \star \rho_+^{-h_2})(t,s)&\le&C\rho_+(t,s)^{- h_1 \wedge h_2}, \hskip1cm  h_i > 1+ \varpi, \, i=1,2.  \label{prop13}
\end{eqnarray}
Moreover, inequalities \eqref{prop11}-\eqref{prop13} are also valid for `discrete' convolution $[\rho_+^{-h_1} \star \rho_+^{-h_2}](t,s),
(t,s) \in \Z^2$ with $\rho(t,s)$ on the r.h.s. of  \eqref{prop11} replaced
by $\rho_+(t,s)$.

\smallskip

\noi (iii) Let $a_i = a_i(t,s), (t,s) \in \Z^2, i=1,2 $ satisfy
$a_i(t,s) = \rho_+(t,s)^{-h_i} \big(L_i(t/\rho_+(t,s)) + o(1)\big),  |t| +|s| \to \infty$, where  $0 < h_i <1+ \varpi, h_1 + h_2 > 1+ \varpi, $
and $L_i(u) \not\equiv 0, u \in [-1,1]$ are bounded piecewise
continuous functions, $i=1,2. $  Let $a_{i\infty} (u,v) :=  \rho(u,v)^{-h_i} L_i(u/\rho(u,v)), (u,v) \in \R^2, i=1,2. $
Then
\begin{eqnarray} \label{prop14}
&[a_1 \star a_2](t,s) \ = \ \rho_+(t,s)^{1+ \varpi - h_1 -  h_2}\big(L_{12}(\frac{t}{\rho_+(t,s)}) + o(1)\big), \quad |t|+|s| \to \infty,
\end{eqnarray}
where
\begin{equation} \label{L12}
L_{12}(z) \ := \  (a_{1\infty}\star a_{2\infty})(z, (1-z^2)^{\varpi/2}) \ = \
\int_{\R^2}
a_{1\infty} (u,v) a_{2\infty} (u+z, v + (1-z^2)^{\varpi/2}) \d u \d v
\end{equation}
is a bounded continuous function on the interval $z \in [-1,1]$.
Moreover, if  $L_1(z) = L_2(z) \ge 0$ then  $L_{12}(z)$  in \eqref{L12} is strictly positive on $[-1,1]$.

\smallskip

\noi (iv) Let  $b(t,s) :=  \rho_+(t,s)^{-h} \big(L(t/\rho_+(t,s)) + o(1)\big),  |t| +|s| \to \infty, (t,s) \in \Z^2, \
b_\infty(t,s) :=  \rho(t,s)^{-h} L(t/\rho(t,s)), $  $ (t,s) \in \R^2 $ where  $0 < h <1+ \varpi$
and $L(u) \ge 0, u \in [-1,1]$ is a
continuous function. Then for any $\gamma >0$
\begin{equation} \label{Bconv}
B_\lambda(\gamma) \ := \ \sum_{  (t_i,s_i) \in K_{[\lambda, \lambda^\gamma]}, i=1,2}
b(t_1-t_2, s_1-s_2)  \
\sim \ {\cal C}(\gamma) \lambda^{2{\cal H}(\gamma)}, \qquad \lambda  \to \infty,
\end{equation}
where
\begin{eqnarray}\label{HC}
&{\cal H}(\gamma) \ := \
\begin{cases}
1 + \varpi - \frac{h}{2}, \\
1+ \gamma - \frac{\gamma h}{2\varpi},   \\
1 + \frac{\gamma}{2} - \frac{h-\varpi}{2}, \\
1 + \gamma - \frac{h}{2},\\
\frac{1}{2} + \gamma  - \frac{\gamma(h-1)}{2\varpi}, 
\end{cases}
\quad {\cal C}(\gamma) \ := \
\begin{cases}
\int_{(0,1]^4}b_\infty (t_1-t_2,s_1-s_2) \d t_1 \d t_2 \d s_1 \d s_2, &{\rm (I)} \\
\int_{(0,1]^2} b_\infty (0,s_1-s_2) \d s_1 \d s_2, &{\rm (II)}\\
\int_{(0,1]^2\times \R} b_\infty (t_1-t_2,s) \d t_1 \d t_2 \d s, &{\rm (III)} \\
\int_{(0,1]^2}b_\infty (t_1-t_2,0) \d t_1 \d t_2, &{\rm (IV)} \\
\int_{\R \times (0,1]^2} b_\infty (t,s_1-s_2) \d t \d s_1 \d s_2, &{\rm (V)}
\end{cases}
\end{eqnarray}
in respective cases {\rm (I)}: $\gamma = \varpi$, {\rm (II)}: $\gamma > \varpi,  h < \varpi$,
{\rm (III)}: $\gamma > \varpi, h > \varpi$, {\rm (IV)}: $\gamma < \varpi,  h < 1 $ and
{\rm (V)}: $\gamma < \varpi,  h > 1$.

\end{proposition}

\section{Covariance structure of subordinated anisotropic RFs}

In this section
from Proposition~\ref{prop1} with $\varpi = \gamma_0$
we obtain the asymptotic form of the covariance function of
$r_X (t,s) := \E X(0,0) X(t,s) $
and the asymptotics of the variance of anisotropic partial sums $S^X_{\lambda, \gamma} $
of subordinated RF $X = A_k(Y)$.

\begin{proposition} \label{appell}  Let RF $X = A_k(Y)$ satisfy assumptions (A1), (A2) and (A3)$_k$.

\smallskip

\noi (i) Let $1\le k < P$. Then
\begin{eqnarray}\label{rhoas2}
r_X(t,s)
&=&\rho(t,s)^{-kp_1} \big(L_X (
t/\rho(t,s)) + o(1)\big), \qquad |t|+|s| \to \infty,
\end{eqnarray}
where $L_X(z):= (a_\infty \star a_\infty)^k(z,  (1-z^2)^{\gamma_0/2}), \, z \in [-1,1] $  is a strictly positive
continuous function and $a_\infty $ is defined in \eqref{V0}.
Moreover, $X(t,s) = Y^{\bullet k}(t,s) + {\cal Z}(t,s)$, where ${\cal Z}(t,s)$ is defined in \eqref{Ykstar} and
\begin{eqnarray}\label{Zcov}
r_{\cal Z}(t,s)
&=&O(\rho(t,s)^{-2q_1}) \ = \  o(\rho(t,s)^{-kp_1}), \qquad |t|+|s| \to \infty.
\end{eqnarray}

\noi (ii) Let $k > P$. Then
\begin{eqnarray}\label{rhoas22}
r_X(t,s)
&=&O\big(\rho(t,s)^{-(kp_1) \wedge (2q_1)}\big), \qquad |t|+|s| \to \infty.
\end{eqnarray}
\end{proposition}

Clearly, \eqref{rhoas2} implies $C_1\rho(t,s)^{-kp_1} \le  r_X(t,s) \le  C_2 \rho(t,s)^{-kp_1}$ for all
$|t|+|s| > C_3 $ and some $0< C_i < \infty, i=1,2,3 $. The last fact together with Proposition~\ref{prop1} (i) implies
the following corollary.

\begin{corollary} \label{XLRD}
Let $X = A_k(Y), \,  1\le k < P$  be the subordinated RF defined
in Proposition \ref{appell} and satisfying the conditions therein.

\smallskip

\noi (i) Let $1\le k < P$. Then $\sum_{(t,s) \in \Z^2} |r_X(t,s)| = \infty $. Moreover,
$\sum_{ s \in \Z} |r_X(0,s)| = \infty \ \Longleftrightarrow \  kp_2 \le 1 $ and
$\sum_{t \in \Z} |r_X(t,0)| = \infty \ \Longleftrightarrow \  kp_1 \le 1. $

\smallskip

\noi (ii) Let $k > P$. Then $\sum_{(t,s) \in \Z^2} |r_X(t,s)| < \infty $.

\end{corollary}

\begin{remark} \label{remLRD}
Following the terminology in \cite{pils2016}, we say that
a covariance stationary RF $X = \{ X(t,s), (t,s) \in \Z^2 \} $ has {\it vertical LRD property}
(respectively, {\it horizontal LRD property}) if
$\sum_{s \in \Z} |r_X(0,s)| = \infty $ (respectively,
$\sum_{t\in \Z} |r_X (t,0)| = \infty $). From
Corollary \ref{XLRD} we see the dichotomy of the limit distribution in Theorems~\ref{thm2} - \ref{thm3}
at points $kp_2 =1 $
at $kp_1 =1 $ is related to the change of vertical and horizontal LRD properties
of the subordinated RF $X = A_k(Y)$.

\end{remark}

\begin{corollary} \label{appellvar} Let $X(t,s) = A_k(Y(t,s)) =   Y^{\bullet k}(t,s) + {\cal Z}(t,s), \,  1\le k < P$  be the subordinated RF defined
in Proposition \ref{appell} and satisfying the conditions therein.
Then for any $\gamma >0$
\begin{eqnarray}\label{varX}
{\rm Var}\big( S^X_{\lambda, \gamma}
\big)\ \sim \ {\rm Var}\big( S^{Y^{\bullet k}}_{\lambda, \gamma}\big)
&\sim&c(\gamma)\lambda^{2H(\gamma)}, \qquad \lambda \to \infty
\end{eqnarray}
and
\begin{equation}
{\rm Var}(S^{\cal Z}_{\lambda,\gamma})
\ =\   O(\lambda^{1+\gamma}), \label{varZ}
\end{equation}
where $H(\gamma) \in ((1+ \gamma)/2, 1+ \gamma)$ and $ c(\gamma)$ are defined in Theorems~\ref{thm1}-\ref{thm3}.

\end{corollary}

\section{Proofs of Theorems \ref{thm1}-\ref{thm5}}

We use the criterion in Proposition~\ref{off} for the convergence in distribution of off-diagonal polygonal forms
towards It\^o-Wiener integral which is a straightforward extension of (\cite{gir2012}, Proposition~14.3.2).

Let $L^2(\Z^{2k}) $ be the class of all real functions $g = g((u,v)_k), (u,v)_k \in \Z^{2k}$ with
$\sum_{(u,v)_k \in \Z^{2k}}  g((u,v)_k)^2 < \infty $ and
$Q_k(g) := \sum^{\bf \bullet}_{(u,v)_k} g((u,v)_k) \vep(u_1,v_1) \cdots \vep(u_k,v_k), g \in L^2(\Z^{2k})
$ be a $k$-tuple
off-diagonal form in i.i.d. r.v.'s $\{\vep(u,v) \} $ satisfying Assumption (A1).
For $g_{\lambda,\gamma} \in L^2(\Z^{2k}) \ (\lambda >0, \gamma >0)$ define a step function
$\widetilde g_{\lambda, \gamma} \in L^2(\R^{2k}) $ by
\begin{eqnarray}\label{tildeg}
\widetilde g_{\lambda,\gamma}((u,v)_k)&:=&\lambda^{k \gamma (1+\gamma_0^{-1})/2}
g_{\lambda,\gamma}([\lambda^{\gamma/\gamma_0}  u_1],[\lambda^\gamma v_1], \cdots,
[\lambda^{\gamma/\gamma_0} u_k],[\lambda^\gamma v_k]), \qquad (u,v)_k \in \R^{2k}.
\end{eqnarray}

\begin{proposition} \label{off}  Assume that there exists  $h_\gamma \in L^2(\R^{2k})$ such that
$ \lim_{\lambda \to \infty} \|\widetilde g_{\lambda, \gamma} - h_\gamma\|_k \to 0.$  Then
$Q_k (g_{\lambda,\gamma}) $   $ \limd \int_{\R^{2k}} h_\gamma((u,v)_k) \d^k W \ (\lambda \to \infty)$.

\end{proposition}

\noi {\it Proof of Theorem \ref{thm1}.}  Let $\rho(t,s) := (|t|^2 + |s|^{2/\gamma_0})^{1/2}, \, (t,s) \in \R^2$.

\smallskip

\noi (i) Let us show that the stochastic integral $V_{k,\gamma_0}(x,y)$  is well-defined
or $\| h(x,y; \cdot)\|_k < \infty $.  where $h(x,y; (u,v)_k)$ is defined in \eqref{V0}.
It suffices to consider the case $x=y = 1$.
By Proposition~\ref{prop1},
$\|h(1,1;\cdot)\|^2_k = \int_{(0,1]^4} ((a_\infty \star a_\infty)(t_1-t_2,s_1 -s_2))^k \d t_1 \d t_2 \d s_1 \d s_2
\le \int_{(0,1]^4} \rho(t_1-t_2,s_1 -s_2)^{-kp_1} \d t_1 \d t_2 \d s_1 \d s_2  < \infty $
since $kp_1 < 1+ \gamma_0 = 1 + p_1/p_2$ or $k < P$ holds.
The self-similarity property
in \eqref{Vss} follows by scaling properties
$a_\infty (\lambda t, \lambda^{\gamma_0} s) = \lambda^{q_1} a_\infty (t,s), \, \{ W(\d \lambda u, \d \lambda^{\gamma_0} v) \}
\eqfdd \{ \lambda^{(1+ \gamma_0)/2} W(\d u, \d v) \}
$ of the integrand
and the white noise, and the change of variables rules for multiple  It\^o-Wiener integral, see
\cite{dobr1979}, also (\cite{gir2012}, Proposition~14.3.5).

\smallskip

\noi (ii) Relation  \eqref{varY} is proved in Proposition~\ref{appellvar}.
Let us prove \eqref{limY}.
Recall the decomposition $X(t,s) = Y^{\bullet k}(t,s) + {\cal Z}(t,s)$ in Corollary~\ref{appellvar}. Using
${\rm Var}(S^{\cal Z}_{\lambda,\gamma_0}) = O(\lambda^{1+\gamma_0}) =  o(\lambda^{2H(\gamma_0)})$,
see \eqref{varZ}, relation \eqref{limY} follows from
\begin{eqnarray}\label{limQ}
Q_k(g_{\lambda,\gamma_0}(x,y; \cdot)) \ = \   \lambda^{-H(\gamma_0)} \sum_{(t,s)\in  K_{[\lambda x, \lambda^{\gamma_0}y]}} Y^{\bullet k}(t,s)
&\limfdd&V_{\gamma_0}(x,y),
\end{eqnarray}
where 
\begin{equation} \label{glambda}
g_{\lambda,\gamma_0}(x,y; (u,v)_k) \ := \  \lambda^{-H(\gamma_0)} \sum_{(t,s)\in  K_{[\lambda x,\lambda^{\gamma_0} y]}}  a(t-u_1,s-v_1) \cdots a(t-u_k,s-v_k),
 \qquad (u,v)_k \in \Z^{2k}.
\end{equation}
Using Proposition~\ref{off} and Cram\'er-Wold device, relation \eqref{limQ} follows  from
\begin{equation}\label{hconv}
\lim_{\lambda \to \infty} \|\sum\nolimits_{i=1}^m \theta_i \big(\widetilde g_{\lambda,\gamma_0}(x_i,y_i; \cdot) - h(x_i,y_i; \cdot)\big)\|_k \ = \ 0,
\end{equation}
for any $m\ge 1$ and any $\theta_i \in \R, (x_i,y_i) \in \R^2_+, 1\le i \le m$,
where the limit function $h(x,y; (u,v)_k)$ is given in \eqref{V0}. We restrict the subsequent proof of
\eqref{hconv} to  the case $m=\theta_1 = 1, (x_1,y_1) = (x,y) $ since the general case
of \eqref{hconv} follows analogously.
Using \eqref{acoefL}, \eqref{glambda}, \eqref{tildeg}
and notation $a_\lambda(t,s) := (\lambda^{-1} \rho(t,s))^{-q_1} \big(L_0(t/\lambda^{-1} \vee \rho(t,s)) + o(1)\big), \, \lambda \to  \infty $
and $ \lambda' := \lambda^{\gamma_0}$ similarly
to \eqref{Glambda} we get
\begin{eqnarray}
\widetilde g_{\lambda,\gamma_0}(x,y; (u,v)_k)&=&\int_{\R^2} \prod_{i=1}^k a_\lambda \big(
\mbox{$\frac{ [\lambda t] - [\lambda u]}{\lambda}, \frac{ [\lambda' s]  -  [\lambda' v]}{\lambda'} $} \big)
\1 \big( ([\lambda t], [\lambda' s]) \in (0, \lambda x] \times (0, \lambda' y]\big) \nn \\
&\to&h(x,y; (u,v)_k) \label{wtig}
\end{eqnarray}
point-wise for any $(u,v)_k \in \R^{2k},  (u_i,v_i) \ne (u_j,v_j) \, (i\ne j)$ fixed.  We use a similar bound to \eqref{a0bdd}, viz.,
\begin{equation}\label{a1bdd}
\mbox{$ \frac{1}{\lambda} \vee \rho \big(\frac{ [\lambda t] - [\lambda u]}{\lambda}, \frac{ [\lambda' s]  -  [\lambda' v]}{\lambda'}\big)$}
\  \ge \ c \rho(t-u, s -v), \quad \forall  \  t,u,s,v \in \R, \quad \exists c >0,
\end{equation}
implying the dominated bound
\begin{eqnarray*}
|\widetilde g_{\lambda,\gamma_0}(x,y; (u,v)_k)|&\le C \int_{(0,2x]\times (0,2y]} \prod_{i=1}^k
\rho(t-u_i,s-v_i)^{-q_1}   \d t \d s \ =: \  \bar g(x,y: (u,v)_k)
\end{eqnarray*}
with $\| \bar g(x,y; \cdot) \|_k < \infty $ so that \eqref{hconv} follows from \eqref{wtig} and Proposition~\ref{prop1}
by the dominated
convergence theorem. Theorem~\ref{thm1} is proved. \hfill $\Box$

\medskip


\noi {\it Proof of Theorem \ref{thm2}.} As noted in Sec.~3, part (iii) follows by the same argument as part (ii)
by exchanging the coordinates $t$ and $s$ and we omit the details.

\smallskip

\noi (i) Let us show that the stochastic integral  in \eqref{Zpm}  is well-defined
or $\| h_+(y; \cdot)\|_k < \infty $, where $h_+(y; (u,v)_k)$ is defined in \eqref{hpm}.
Indeed by Proposition~\ref{prop1}
$\|h_+ (y; \cdot)\|^2_k = \int_{(0,1]^2} ((a_\infty \star a_\infty)(0,s_1 -s_2))^k  \d s_1 \d s_2
\le\ C\int_{(0,1]^2} \rho(0,s_1 -s_2)^{-kp_1}  \d s_1 \d s_2
\le C\int_{[-1,1]} |s|^{-kp_2}  \d s \ < \  \infty $
since $kp_2 < 1 $. The remaining facts in (i) follow
similarly as in the proof of Theorem \ref{thm1}(i).

\smallskip

\noi (ii)  Relation  \eqref{varY2} is proved in Corollary~\ref{appellvar}.
Similarly to the proof of \eqref{limY}, the weak convergence in \eqref{limY2} follows from
\begin{eqnarray}\label{limQ2}
Q_k(g_{\lambda,\gamma}(x,y; \cdot)) \ = \   \lambda^{-H(\gamma)} \sum_{(t,s)\in  K_{[\lambda x, \lambda^{\gamma}y]}} Y^{\bullet k}(t,s)
&\limfdd&x Z_+(y),
\end{eqnarray}
where 
\begin{equation} \label{glambda2}
g_{\lambda,\gamma}(x,y; (u,v)_k) \ := \  \lambda^{-H(\gamma)} \sum_{(t,s)\in  K_{[\lambda x,\lambda^{\gamma} y]}}
a(t-u_1,s-v_1) \cdots a(t-u_k,s-v_k),
 \qquad (u,v)_k \in \Z^{2k}.
\end{equation}
Again, we restrict the proof of \eqref{limQ2} to one-dimensional convergence at $(x,y)\in \R^2_+$.
By Proposition \ref{off} this follows from
\begin{equation}\label{hconv2}
\lim_{\lambda \to \infty} \|\widetilde g_{\lambda,\gamma}(x,y; \cdot) - x h_+(y; \cdot)\|_k \ = \ 0,
\end{equation}
where, with $\lambda' := \lambda^\gamma, \lambda'' := \lambda^{\gamma/\gamma_0} \gg \lambda,
a_{\lambda''} (t,s) := ( (\lambda'')^{-1} \vee \rho (t,s) )^{-q_1} ( L_0 (t/ ( (\lambda'')^{-1} \vee \rho(t,s) ) ) + o(1) ) $,
\begin{eqnarray}
\widetilde g_{\lambda,\gamma}(x,y; (u,v)_k)&=&\int_{\R^2} \prod_{i=1}^k
a_{\lambda''} \big( \textstyle\frac{[\lambda t] - [\lambda'' u_i]}{\lambda''}, \frac{[\lambda' s]  -  [\lambda' v_i]}{\lambda'} \big)
\1 \big( ([\lambda t], [\lambda' s]) \in (0, \lambda x] \times (0, \lambda' y]\big) \d t \d s \nn \\
&\to&x  h_+(y; (u,v)_k) \label{wtig2}
\end{eqnarray}
point-wise for any $(u,v)_k \in \R^{2k},  (u_i,v_i) \ne (u_j,v_j) \, (i\ne j)$ fixed.

The  dominating convergence argument to
prove  \eqref{hconv2} from \eqref{wtig2} uses Pratt's  lemma \cite{pra1960}, as follows. Similarly to  \eqref{a1bdd}
note that
\begin{equation*}\label{a2bdd}
\mbox{$ \frac{1}{\lambda''} \vee \rho \big(\frac{ [\lambda t] - [\lambda'' u]}{\lambda''}, \frac{ [\lambda' s]  -  [\lambda' v]}{\lambda'}\big)$}
\  \ge \ c \rho((\lambda t/\lambda'') -u, s -v), \quad \forall  \  t,u,s,v \in \R, \quad \exists c >0
\end{equation*}
and hence
\begin{eqnarray*}
|\widetilde g_{\lambda,\gamma}(x,y; (u,v)_k)|&\le&C\int_{(0,2x]\times (0,2y]} \prod_{i=1}^k
 \rho( (\lambda t/\lambda'') - u_i, s  - v_i)^{-q_1} \d t \d s
   \ =: \ C G_\lambda ((u,v)_k)
\end{eqnarray*}
with $C>0 $ independent of $\lambda >0, (u,v)_k \in \R^{2k}$. Clearly,
$\lim_{\lambda \to  \infty} G_\lambda ((u,v)_k) = G ((u,v)_k) := 2x \int_{(0,2y]} \prod_{i=1}^k
$ $\rho(- u_i, s  - v_i)^{-q_1}\d s $ point-wise in $\R^{2k} $ and
\begin{eqnarray*}
 \|  G_\lambda \|^2_k
 &=&\int_{(0,2x]^2 \times (0,2y]^2} \big( (\rho^{-q_1} \star \rho^{-q_1})( (\lambda/\lambda'') (t_1-t_2), s_1-s_2) \big)^k
 \d t_1 \d t_2 \d s_1 \d s_2 \\
&\to&\int_{(0,2x]^2 \times (0,2y]^2} \big( (\rho^{-q_1} \star \rho^{-q_1})(0, s_1-s_2) \big)^k
 \d t_1 \d t_2 \d s_1 \d s_2 \ = \  \| G\|^2_k < \infty
 \end{eqnarray*}
by Proposition \ref{prop1} and condition $1\le k < 1/p_2$, or $ p_2 = q_2(2-Q) < 1/k$. Thus,
application of \cite{pra1960} proves  \eqref{hconv2}.  Theorem \ref{thm2}  is proved.
\hfill $\Box$

\medskip

To prove Theorem \ref{thm3} we use approximation by $m$-dependent variables and the following CLT for triangular array
of $m$-dependent r.v.'s.

\begin{lem} \label{cltm} Let $\{\xi_{ni}, 1\le i \le N_n \}, n \ge 1 $ be a triangular array of $m$-dependent r.v.'s with zero mean
and finite variance.
Assume that:  (L1) $\xi_{in}, 1 \le i \le N_n $ are identically distributed for any $n\ge 1$,
(L2) $ \xi_{n}  :=  \xi_{n1}  \limd  \xi,  \, \E \xi^2_n \to \E \xi^2 < \infty  $ for some r.v. $\xi $ and
(L3) ${\rm Var} (\sum_{i=1}^{N_n} \xi_{ni}) \sim \sigma^2 N_n, \, \sigma^2 >0$.
Then $N_n^{-1/2} \sum_{i=1}^{N_n} \xi_{ni} \limd N(0, \sigma^2)$.
 \end{lem}

\noi {\it Proof.} W.l.g., we can assume $N_n = n$ in the subsequent proof.
We use the CLT due to Orey \cite{orey1958}. Accordingly, let
$\xi^\tau_{ni} := \xi_{ni} \1( |\xi_{ni}| \le \tau n^{1/2}), \alpha^\tau_{ni} := \E \xi^\tau_{ni},
\sigma^\tau_{nij} := {\rm Cov}(\xi^\tau_{ni}, \xi^\tau_{nj})$. It suffices to show that for any $\tau >0$ the following
conditions in \cite{orey1958} are satisfied:
(O1) $n^{-1/2}\sum_{i=1}^{n} \alpha^\tau_{ni} \to 0$, (O2) $n^{-1}\sum_{i,j=1}^{n} \sigma^\tau_{nij} \to \sigma^2, $
(O3) $n^{-1}\sum_{i,j=1}^{n} \sigma^\tau_{nii} = O(1), $ and (O4)
$\sum_{i=1}^{n} \P (|\xi_{ni}| >  \tau n^{1/2}) \to 0. $ \\
Consider (O1), or $n^{1/2} \alpha^\tau_n \to 0, \alpha^\tau_n := \alpha^\tau_{n1}.  $ We have
$0 =  n^{1/2} \E \xi_n =  n^{1/2} \alpha^\tau_n + \kappa_n, $ where
$|\kappa_n| := n^{1/2} | \E \xi_{n} \1( |\xi_{n}| > \tau n^{1/2}) | \le  \tau^{-1} \E \xi^2_n \1 ( |\xi_n| > \tau n^{1/2}). $
Therefore, (O1) follows from
\begin{equation} \label{xi2}
\E \xi^2_n \1 ( |\xi_n| > \tau n^{1/2}) \ \to \ 0.
\end{equation}
Using the Skorohod representation theorem \cite{sko1956} w.l.g. we can assume that r.v.s
$\xi, \xi_n, n \ge 1 $ are defined on the same probability space and
$\xi_n \to \xi $ almost surely. The latter fact together with (L2) and Pratt's  lemma \cite{pra1960}
implies that $\E |\xi^2_n - \xi^2| \to 0$ and hence \eqref{xi2} follows
due to $ \P ( |\xi_n| > \tau n^{1/2}) \to 0,$  see (\cite{nev1964}, Ch.2, Prop.5.3). The above
argument also implies (O4) since $\P (|\xi_{n}| >  \tau n^{1/2}) \le \tau^{-1} n^{-1} \E \xi^2_n \1 ( |\xi_n| > \tau n^{1/2})$
by Markov's inequality. (O3) is immediate from (L1) and (L2). Finally,
(O2) follows from (L3), (O1) and
\begin{equation} \label{xi3}
n^{-1} \sum_{1\le i,j\le n, |i-j| \le m}
\E (\xi_{ni} \xi_{nj} - \xi^\tau_{ni} \xi^\tau_{nj}) \ \to \ 0.
\end{equation}
Let $ \widetilde \xi^\tau_{ni} := \xi_{ni} - \xi^\tau_{ni}$.
Since $|\E (\xi_{ni} \xi_{nj} - \xi^\tau_{ni} \xi^\tau_{nj})| \le
| \E (\widetilde \xi^\tau_{ni} \xi^\tau_{nj} +  \xi^\tau_{ni} \widetilde \xi^\tau_{nj} + \widetilde \xi^\tau_{ni} \widetilde \xi^\tau_{nj})|
\le C \E^{1/2} \xi^2_n \1 ( |\xi_n| > \tau n^{1/2}) $, relation \eqref{xi3} follows from \eqref{xi2}.
Lemma \ref{cltm} is proved. \hfill $\Box$

\vskip.3cm

\noi {\it Proof of Theorem \ref{thm3}.} Again, we prove part (i) only since part (ii)
follows similarly by exchanging the coordinates $t$ and $s$.  \\
Relation  \eqref{varY3} is proved in Proposition \ref{appellvar}.  Let us prove
\eqref{limY3}. Similarly as in the case of the previous theorems, we shall restrict ourselves
with the proof of one-dimensional convergence at $(x,y) \in \R^2_+$.
For
$m \ge 1, \lambda >0 $ define stationary RFs
\begin{eqnarray}
X_{m}(t,s)&:=&A_k(Y_{m}(t,s)), \qquad \text{where} \nn \\
Y_{m}(t,s)&:=&
\sum_{(u,v) \in \Z^2: |s-v| \le [\lambda^{\gamma_0}]m} a(t-u,s-v) \vep(u,v),   \label{Ym1}
\end{eqnarray}
and where $A_k $ stands for the Appell polynomial of degree $k$ relative to the distribution of $Y_{m}(t,s) $.
Note $X_{m}(t_1,s_1) $ and $X_{m}(t_2,s_2) $ are independent if $|s_1 - s_2| > 2 [\lambda^{\gamma_0}] m $.
Then
\begin{eqnarray}\label{Usum}
&S^{X_m}_{\lambda,\gamma}(x,y)\ := \ \sum_{(t,s)\in  K_{[\lambda x, \lambda^{\gamma}y]}} X_{m}(t,s) \ = \ \sum_{i=0}^{N} U_{\lambda,m}(i)
\end{eqnarray}
where $N := [[\lambda^\gamma y]/[\lambda^{\gamma_0}]] = O(\lambda^{\gamma - \gamma_0}) $ and
\begin{eqnarray}\label{Snm}
U_{\lambda,m}(i)
&:=&\sum_{1\le t \le \lambda x} \ \sum_{  i [\lambda^{\gamma_0}] < s \le (i+1) [\lambda^{\gamma_0}] } X_{m}(t,s).  \label{Unm}
\end{eqnarray}
Note $U_{\lambda,m}(i)$ and $U_{\lambda,m}(j)$ are independent provided $|i-j| > 2m$ hence \eqref{Usum} is a sum
of $2m$-dependent r.v.'s.
The one-dimensional convergence in \eqref{limY3} follows from standard Slutsky's argument (see e.g. \cite{gir2012}, Lemma 4.2.1) and
the following lemma. Theorem~\ref{thm3} is proved. \hfill $\Box$

\begin{lem} \label{lemSm}
Under the conditions and notation of Theorem \ref{thm3} (i), for any $\gamma > \gamma_0$ and any $m = 1,2, \cdots $
\begin{eqnarray}\label{S1}
{\rm Var}(S^{X_m}_{\lambda,\gamma}(x,y))&\sim&\sigma^2_m(x,y) 
\lambda^{2H(\gamma)} \qquad \text{and} \qquad \\
\lambda^{-H(\gamma)} S^{X_m}_{\lambda,\gamma}(x,y)&\limd&N(0, \sigma^2_m(x,y)) \quad \text{as} \quad \lambda \to \infty,
\label{S2}
\end{eqnarray}
where $\sigma^2_m(x,y) $ is defined in \eqref{S1m} below.  Moreover,
\begin{equation}
\lim_{m\to \infty}  \limsup_{\lambda \to \infty} \lambda^{-2H(\gamma)} {\rm Var}(S^X_{\lambda,\gamma}(x,y) - S^{X_m}_{\lambda,\gamma}(x,y)) \ = \  0. \label{S3}
\end{equation}
\end{lem}

\noi {\it Proof.} By adapting the argument in the proof of \eqref{varY3} and
Proposition \ref{prop1} (iv), Case (III), we can show the limits
\begin{eqnarray}\label{S1m}
\lambda^{-2H(\gamma)} {\rm Var}(S^{X_m}_{\lambda,\gamma}(x,y))
&\to&y\int_{(0,x]^2 \times \R}
\big((a_{\infty,m} \star a_{\infty,m})(t_1-t_2,s)\big)^k  \d t_1 \d t_2  \d s \ =: \ \sigma^2_m(x,y)
\end{eqnarray}
and
\begin{eqnarray}
&&\lambda^{-2H(\gamma)} {\rm Var}(S^X_{\lambda,\gamma}(x,y) - S^{X_m}_{\lambda,\gamma}(x,y)) \nn \\
&&= \ \lambda^{-2H(\gamma)}\sum_{(t_i,s_i) \in K_{[\lambda x, \lambda^\gamma y]}, i=1,2}
\big\{{\rm Cov}(X(t_1,s_1), X(t_2,s_2)) -  {\rm Cov}(X(t_1,s_1), X_m(t_2,s_2)) \nn\\
&& \hskip5cm  - \, {\rm Cov}(X_m(t_1,s_1), X(t_2,s_2)) + {\rm Cov}(X_m(t_1,s_1), X_m(t_2,s_2))\big\}\nn  \\
&&\to\  y\int_{(0,x^2]\times \R} G_m (t_1-t_2,s) \d t_1 \d t_2 \d s, \qquad  \lambda \to \infty,   \label{S3m}
\end{eqnarray}
where  $G_m(t,s):= ((a_\infty \star a_\infty)(t,s))^k -  ((a_\infty \star a_{\infty,m})(t,s))^k  -
((a_{\infty,m} \star a_\infty)(t,s))^k + ((a_{\infty,m} \star a_{\infty,m})(t,s))^k $ and
\begin{equation}\label{aminfty}
a_{\infty,m}(t,s)\ := \ L_0(t/\rho(t,s))\rho(t,s)^{-q_1} \1(|s|\le m), \quad (t,s) \in \R^2.
\end{equation}
is a `truncated' version of $a_\infty (t,s)$ in \eqref{V0}. Since $|G_m(t,s)| \le 4  (a_\infty \star a_\infty)(t,s))^k $ and
$G_m(t,s)$ vanishes with $m \to \infty $ for any fixed $(t,s) \ne (0,0) $, \eqref{S3}  follows from
\eqref{S3m} by the dominated convergence theorem.

The proof of  \eqref{S2} uses Lemma \ref{cltm}. Accordingly,
let $N_\lambda := [[\lambda^\gamma y]/[\lambda^{\gamma_0}]] $ and $\xi_{\lambda i} := \lambda^{-H(\gamma_0)} U_{\lambda,m}(i)$, where $H(\gamma_0) = 1 + \gamma_0 - kp_1/2 $ is
the same as in Thm \ref{thm1} and
$U_{\lambda,m}(i) $ are $2m$-dependent r.v.'s defined in \eqref{Unm}. Note $U_{\lambda,m}(i), i=1, \cdots, N_\lambda -1$ are identically distributed and
$\lambda^{H(\gamma_0)}N^{1/2}_\lambda \sim \lambda^{H(\gamma)} y^{1/2}$. Thus, condition (L1)
of Lemma \ref{cltm} for $\xi_{\lambda i}, 1 \le i \le N_\lambda -1 $
is satisfied and
(L3) follows from
${\rm Var}(\sum_{i=1}^{N_\lambda -1} \xi_{\lambda i}) \sim  {\rm Var} (S^X_{\lambda,\gamma}(x,y))  \sim
c_m(\gamma) x^{2H_{1k}} y,  $
see \eqref{S1}. Finally, condition (L2), or
\begin{equation}
\xi_{\lambda,1}  \ = \  \lambda^{-H(\gamma_0)} U_{\lambda,m}(1) \ \limd \ \xi, \qquad \E \xi^2_{\lambda,1} \ \to \ \E \xi^2
\end{equation}
follows similarly as in Theorem \ref{thm1} with the limit r.v. $\xi$ given by the $k$-tuple It\^o-Wiener integral:
$$
\xi :=
\int_{\R^{2k}} \Big\{\int_0^x \int_0^1 \prod_{\ell =1}^k a_{\infty,m}(t-u_\ell, s-v_\ell)\, \d t \d s \Big\} \d^k W
$$
and $a_{\infty,m}(t,s)$ defined in  \eqref{aminfty}. This proves  \eqref{S2} and Lemma \ref{lemSm}, too.
\hfill $\Box$

\medskip

\noi {\it Proof of Theorem \ref{thm4}.} The proof is an adaptation of
the proof of CLT in (\cite{gir2012},
Theorem~4.8.1) for sums of  `off-diagonal' polynomial forms with one-dimensional `time' parameter.
Define
\begin{eqnarray}
X_{m}(t,s)&:=&A_k(Y_{m}(t,s)), \qquad   
Y_{m}(t,s)\ := \
\sum_{(u,v) \in \Z^2: |t-u|+|s-v| \le m} a(t-u,s-v) \vep(u,v), \label{Ym2}
\end{eqnarray}
where $A_k $ stands for the Appell polynomial of degree $k$ relative to the distribution of $Y_{m}(t,s) $.
Note the truncation level  $m$ in \eqref{Ym2} does not depend on $\lambda $ in contrast to  the truncation level
$m [\lambda^{\gamma_0}]$ in \eqref{Ym1}. Similarly to Lemma \ref{lemSm} it suffices to prove for
any $\gamma >0, m = 1,2, \cdots $
\begin{eqnarray}
&{\rm Var}(S^{X_m}_{\lambda,\gamma}(x,y))\ \sim\  x y\sigma^2_{X_m} 
\lambda^{1+ \gamma}, \quad 
\lambda^{-(1+\gamma)/2} S^{X_m}_{\lambda,\gamma}(x,y)
\ \limd\ N(0, x y \sigma^2_{X_m}) 
\label{S4} \\
&\lim_{m\to \infty}  \limsup_{\lambda \to \infty} \lambda^{-(1+\gamma)}
{\rm Var}(S^X_{\lambda,\gamma}(x,y) - S^{X_m}_{\lambda,\gamma}(x,y)) \ = \  0. \label{S5}
\end{eqnarray}
where $\sigma^2_{X_m} := \sum_{(t,s) \in \Z^2} r_{X_m} (t,s) $ and
$r_{X_m}(t,s) := {\rm Cov}(X_m(0,0), X_m(t,s))$. Note $X_{m}(t_1,s_1)$ and $X_{m}(t_2,s_2)$ are independent if
$|t_1-t_2| + |s_1 - s_2| > 2m $. Therefore $\sum_{(t,s) \in \Z^2} |r_{X_m} (t,s)| < \infty $ and
\eqref{S4}  follows the CLT for $m$-dependent RFs, see  \cite{bolt1982}
Consider \eqref{S5}, where we can put $x=y=1$ w.l.g. We have
$\lambda^{-(1+ \gamma)}{\rm Var}(S^X_{\lambda,\gamma} - S^{X_m}_{\lambda,\gamma}) \le $  $ \sum_{(t,s) \in \Z^2} |\phi_m(t,s)|,  $ where
$\phi_m(t,s) := {\rm Cov}\big(X(0,0) - X_m(0,0), X(t,s) - X_m(t,s)). $
From \eqref{Zap1}, \eqref{Zbdd} and \eqref{Ykcov}
we conclude that
$$
|{\rm Cov} (X(0,0), X(t,s))| +
|{\rm Cov} (X(0,0), X_m(t,s))| +  |{\rm Cov} (X_m(0,0), X_m(t,s))|
\le C \rho(t,s)^{-(kp_1) \wedge (2q_1)}
$$
as in \eqref{rhoas22}, with $C >0 $ independent of $m$. Therefore,
$ |\phi_m(t,s)| \le C  \rho(t,s)^{-(kp_1) \wedge (2q_1)} =: \phi(t,s)$, where
$\sum_{(t,s) \in \Z^2} \phi(t,s) < \infty $, see Proposition \ref{prop1}(i), also Corollary \ref{XLRD}(ii).
Thus,  \eqref{S5} follows by the dominated convergence theorem and the fact that
$\lim_{m \to \infty} \phi_m(t,s) = 0$ for any $(t,s) \in \Z^2 $.
Theorem \ref{thm4} is proved.
\hfill $\Box$

\medskip

\noi {\it Proof of Theorem \ref{thm5}.} (i) Split $X = X_k  + X'_k,$ where $X'_k := \sum_{j=k+1} c_j X_j/j!, \,
X_j(t,s) := H_j(Y(t,s))$. Since all statements
of Theorems \ref{thm1}-\ref{thm3}  hold for RF $X_k = H_k(Y)$ and ${\rm Cov}(X_k(t_1,s_1),  X'_k(t_2,s_2)) = 0, \, \forall
(t_i,s_i) \in \Z^2, i=1,2,$
it suffices to show that
\begin{equation} \label{Xrem}
{\rm Var}(S^{X'_k}_{\lambda,\gamma}) \ = \ o(\lambda^{2H(\gamma)}), \qquad  \lambda \to \infty
\end{equation}
for $H(\gamma)$ defined in  Theorems \ref{thm1}-\ref{thm3}.  By well-known properties of Hermite polynomials,
${\rm Var}(S^{X'_k}_{\lambda,\gamma}) = \sum_{j=k+1}^\infty c^2_j {\rm Var}(S^{X_j}_{\lambda,\gamma})/(j!)^2 $
and ${\rm Var}(S^{X_j}_{\lambda,\gamma}) =  j!\sum_{(t_i,s_i) \in K_{[\lambda,\lambda^\gamma]}, i=1,2}
r^j_Y (t_1-t_2, s_1-s_2) \le j! \Sigma_{k+1}(\lambda)$, where
$\Sigma_{k+1}(\lambda) $  $ := \sum_{(t_i,s_i) \in K_{[\lambda,\lambda^\gamma]}, i=1,2}
|r_Y (t_1-t_2, s_1-s_2)|^{k+1}
$ for $j \ge k+1 $ since
$|r_Y(t,s)| \le 1$ according to Assumption (A4)$_k$. Hence, ${\rm Var}(S^{X'_k}_{\lambda,\gamma}) \le \big(\sum_{j=k+1}^\infty c^2_j/j!\big)
\Sigma_{k+1}(\lambda) \le \E G(Y(0,0))^2 \Sigma_{k+1}(\lambda) $, where
$\Sigma_{k+1}(\lambda) = o(\lambda^{2H(\gamma)})$ follows by Proposition \ref{prop1}. This proves
\eqref{Xrem} and part (i).

\smallskip

\noi (ii) For large $ K \in \N,  K >k$, split $ X =  X_K + X'_K$, where $X_K(t,s) =  \sum_{j=k}^K c_j H_j(Y(t,s))/j! $ and
$X'_K(t,s) :=  \sum_{j=K+1}^\infty c_j H_j(Y(t,s))/j!.$  Then ${\rm Var}(S^{X'_K}_{\lambda,\gamma})
\le \big(\sum_{j=K+1}^\infty c^2_j/j!\big)
\Sigma_{K+1}(\lambda)$ as in the proof of part (i), implying  ${\rm Var}(S^{X'_K}_{\lambda,\gamma}) \le C \epsilon_K \lambda^{1 + \lambda } $
where $\epsilon_K := \sum_{j=K+1}^\infty c^2_j/j! $ can be made arbitrary small by choosing $K$ large enough.
On the other hand,  by Theorem \ref{thm4} $\lambda^{-(1+\gamma)/2} S^{X_j}_{\lambda,\gamma}(x,y) \limfdd \sigma_{X_j} B_{1/2,1/2}(x,y)$ for
any $j \ge k$ and the last result extends to finite sums of Hermite polynomials, viz.,
$\lambda^{-(1+\gamma)/2} S^{X_K}_{\lambda,\gamma}(x,y) \limfdd \sigma_{X_K} B_{1/2,1/2}(x,y)$, where
$\sigma^2_{X_K} =  \sum_{(t,s)\in \Z^2} {\rm Cov}(X_K(0,0), X_K(t,s)) \to \sigma^2_X, \, K \to \infty. $
See e.g. (\cite{gir2012}, proof of Thm.4.6.1). The remaining details are easy.
Theorem \ref{thm4} is proved. \hfill $\Box$

\section{Proofs of Propositions \ref{propex2}, \ref{prop1}, \ref{appell} and Corollary \ref{appellvar} }

{\it Proof of Proposition \ref{propex2}.}
The transition probabilities $q_u(v)$ in \eqref{auv} can be explicitly written in terms of binomial probabilities
${\rm bin}(j, k; p) := {k \choose j} p^j (1-p)^{k-j}, \, k = 0,1,\cdots, j =0,1,\cdots,k, \, 0 \le p \le 1$:
\begin{equation}\label{quv}
q_u(v) \ = \  \sum_{j=0}^u {\rm bin}(u-j, u; \theta) {\rm bin} ( (v+j)/2, j; 1/2), \qquad u \in \N, \ |v| \le u.
\end{equation}
Similarly to (\cite{kou2016}, proof of Prop.4.1)
we shall use the following version of the Moivre-Laplace theorem (Feller \cite{fell1966},  ch.7, \S 2, Thm.1):
{\it There exists a constant $C$  such
when $j \to \infty $ and $k \to \infty$ vary in such a way that
\begin{equation} \label{moivre1}
\frac{(j - kp)^3 }{k^2} \ \to \ 0,
\end{equation}
then}
\begin{eqnarray}\label{moivre2}
\bigg|  \frac{{\rm bin}(j,k;p)}{ \frac{1}{\sqrt{2 \pi k p(1-p)}}  \exp \{- \frac{(j -kp)^2}{2k p(1-p)}\}} - 1\bigg|&<&
\frac{C}{k} + \frac{C |j-kp|^3 }{k^2}.
\end{eqnarray}
Let us first explain the idea of the proof.
Using \eqref{quv} and replacing the binomial probabilities by Gaussian densities according to \eqref{moivre2} leads to
\begin{eqnarray*}
a(u,v)
&\sim&\frac{1}{2}\sum_{j=0}^u \frac{1}{\Gamma(d) u^{1-d}}
\frac{1}{\sqrt{2\pi \theta (1-\theta) u}} \e^{- (j- (1-\theta)u)^2/2\theta (1-\theta)u} \,\frac{1}{\sqrt{j \pi/2}} \e^{- v^2/2j} \\
&=&\frac{1}{\Gamma(d) \sqrt{2\pi} u^{3/2 -d}} \sum_{ j=0}^v u^{-1} \frac{1}{\sqrt{2\pi \theta(1-\theta)/u}}
\e^{-  u((j/u)- (1-\theta))^2/2\theta (1-\theta)}  \frac{1}{\sqrt{j/u}}  \e^{- (v^2/u)/2(j/u)} \\
&\sim&\frac{1}{\Gamma(d) \sqrt{2\pi} u^{3/2 -d}}\int_0^1 \frac{1}{\sqrt{2\pi \theta(1-\theta)/u}}
\e^{-  u(x- (1-\theta))^2/2\theta (1-\theta)}  \frac{1}{\sqrt{x}}  \e^{- (v^2/u)/2x} \d x \\
&\sim&\frac{1}{\Gamma(d) \sqrt{2\pi} u^{3/2 -d}} \frac{1}{\sqrt{1-\theta}}  \e^{- (v^2/u)/2(1-\theta)} \\
&=&\rho(u,v)^{d-3/2} L_0(z)
\end{eqnarray*}
with $L_0(z)$ defined in \eqref{L0}.  Here, factor $1/2$ in front of the sum in the first line appears since
${\rm bin} ( (v+j)/2, j; 1/2) = 0$ whenever $v+j$ is odd, in other words, by using Gaussian approximation for all (even and odd)
$j $ we double
the sum and therefore must divide it by $2$. Note also that in the third line, the Gaussian kernel
$ \frac{1}{\sqrt{2\pi \theta(1-\theta)/u}}
\e^{-  u(x- (1-\theta))^2/2\theta (1-\theta)} $ acts as a $\delta$-function at $x = 1-\theta $ when $u \to \infty $.

Let us turn to a rigorous proof of the above asymptotics. For $(u,v) \in \Z^2, (u,v) \ne (0,0),$
denote $\varrho := (u^2 + v^4)^{1/2},
z := u/\varrho  \in [-1,1]$, then $ u = z \varrho, v^2 = \varrho \sqrt{1 - z^2}$.
It suffices to prove
\begin{eqnarray}\label{J0lim}
\varrho^{3/2 -d}a(u,v) - L_0(z)&\to&0 \qquad \text{as}  \  |u|+|v| \to \infty.
\end{eqnarray}
By definition (see \eqref{auv}, \eqref{L0}), \eqref{J0lim}  holds for $u \le 0, z \le 0$ hence
we can assume $u \ge 1, z >0 $ in as follows. Moreover, for any $\epsilon >0$ there exists $ K >0$ such that
\begin{equation}\label{J00}
\varrho^{3/2 -d} a(u,v)  < \epsilon \quad \text{and} \quad L_0(z) < \epsilon \quad (\forall \, 1\le u < v^{9/5}, \ \varrho > K).
\end{equation}
The second relation in \eqref{J00}  is immediate by $\lim_{z \to 0} L_0(z) = L_0(0) =  0$ and $z = u/\varrho \le
\varrho^{9/10}/\varrho \to 0 \ (\varrho
\to \infty)$.
To prove the first relation
we use Hoeffding's inequality  \cite{hoef1963}. Let ${\rm bin}(j, k; p)$ be the binomial distribution. Then
for any $\tau >0$
\begin{equation}\label{hoef}
\sum_{0\le j \le k: |j-kp| > \tau \sqrt{k}}   {\rm bin}(j, k; p) \ \le \ 2 \e^{-2\tau^2}.
\end{equation}
\eqref{hoef} implies ${\rm bin} ( (v+j)/2, j; 1/2) \le  2\e^{-2 v^2/j} \le 2 \e^{-2v^2/u} $ for any $|v| \le u, 0 \le j \le u $.
Also note that $1 \le u \le v^{9/5} $ implies $v^2 \ge 2^{1/20} u \varrho^{1/10} $.
Using  these facts and \eqref{quv} with $\sum_{j=0}^u {\rm bin}(u-j, u; \theta) = 1 $ for any
$1\le u < v^{8/5}$
we obtain
\begin{equation*}
\varrho^{3/2 -d} a(u,v)\ \le \ C\varrho^{3/2 -d} q_u(v) 
\ \le  \ C \varrho^{3/2 -d} \e^{-2v^2/u} \ \le \ C \varrho^{3/2 -d}  \e^{-\varrho^{1/10}} \ \to  \ 0, \quad  \varrho \to \infty,
\end{equation*}
proving  \eqref{J00}.
Hence, it suffices to prove
\eqref{J0lim} for $u \to \infty, 0\le v \le u^{5/9}$.
Below, we give the proof  for $v $ even, the proof for $v$ odd being similar.
Denote
\begin{eqnarray*}
&{\cal D}^+(u,v)\ :=\ \{0 \le j \le u/2:   |2j -  u(1-\theta)| < u^{3/5} \ \text{and} \ |v| < j^{3/5} \}, \\
&{\cal D}^-(u,v)\ :=\ \{0\le j \le u/2:   |2j -  u(1-\theta)| \ge  u^{3/5}\ \ \text{or} \ \ |v| \ge j^{3/5} \}.
\end{eqnarray*}
Split $a(u,v) = \psi_u(-d)\sum_{0\le j \le u/2}  {\rm bin}(u-2j, u; \theta) {\rm bin} (v/2 +j, 2j; 1/2)
= a^+(u,v) + a^-(u,v), $
where  $a^\pm(u,v): = \psi_u(-d) \sum_{j \in {\cal D}^\pm(u,v) }\cdots$.  It suffices to prove that
\begin{eqnarray}\label{Jlim}
\varrho^{3/2 -d}a^+(u,v) - L_0(z)&\to&0 \qquad \text{and} \qquad
\varrho^{3/2 -d}a^-(u,v)\ \to \ 0
\end{eqnarray}
as $u \to \infty, 0\le v \le u^{5/9}$.
To show the first relation in \eqref{Jlim}, let $j^*_u := [u(1-\theta)/2] $ and
\begin{eqnarray*}
a^*(u,v)
&:=&{\rm bin} (v/2 + j^*_u, 2j^*_u; 1/2) \psi_u(-d)
\sum_{j\in {\cal D}^+(u,v)}   {\rm bin}(u-2j, u; \theta),
\end{eqnarray*}
then
\begin{eqnarray*}
a^*(u,v)-a^+(u,v)
&=& \psi_u(-d)
\sum_{j\in {\cal D}^+(u,v)}  {\rm bin}(u-2j, u; \theta)
\big({\rm bin} (v/2 + j^*_u, 2j^*_u; 1/2)- {\rm bin} (v/2 +j, 2j; 1/2)\big).
\end{eqnarray*}
According to \eqref{moivre2},  for $j\in {\cal D}^+(u,v), j^*_u\in {\cal D}^+(u,v)$
\begin{eqnarray*}
{\rm bin} (v/2 +j, 2j; 1/2)&=&\frac{1}{\sqrt{\pi j}} \e^{-v^2/4j}\big(1 + O(j^{-1/5})\big) \ = \
\frac{1}{\sqrt{\pi j}} \e^{-v^2/4j}\big(1 + O(u^{-1/5})\big), \\
{\rm bin} (v/2 + j^*_u, 2j^*_u; 1/2)
&=&\frac{1}{\sqrt{\pi j^*_u}} \e^{-v^2/4j^*_u  }\big(1 + O(u^{-1/5})\big).
\end{eqnarray*}
Using $ c_- u < j < c_+ u,  j\in {\cal D}^+(u,v)$ for some $c_\pm > 0$, and elementary inequalities we obtain that
$|\frac{1}{\sqrt{\pi j}} \e^{-v^2/4j} - \frac{1}{\sqrt{\pi j^*_u}} \e^{-v^2/4j^*_u  }| \le C u^{-7/10} \e^{- c v^2/u } $
for some $C, c >0$ and hence the bound
\begin{eqnarray*}
\big|{\rm bin} (v/2 + j^*_u, 2j^*_u; 1/2)- {\rm bin} (v/2 +j, 2j; 1/2)\big|
&\le&C u^{-7/10} \e^{- c v^2/u }
\end{eqnarray*}
for all $j\in {\cal D}^+(u,v)$ and all $u >0$ large enough. Therefore  since $\sum_{j\in {\cal D}^+(u,v)} {\rm bin}(u-2j, u; \theta) \le 1 $
we obtain
\begin{eqnarray*}
\varrho^{3/2 -d} |a^*(u,v)-a^+(u,v)|
&\le&C\varrho^{3/2 -d} u^{-7/10 + d-1} \e^{- c v^2/u }  \ = \  \varrho^{-1/5} L^*(z) \ \le  \  C \varrho^{-1/5},
\end{eqnarray*}
where $L^*(z)  :=  C z^{d  - 17/10} \e^{-c \sqrt{(1/z)^2 -1}}, z \in (0,1] $ is a bounded function.  As a consequence,
it suffices to prove the first relation in \eqref{Jlim} with $a^+(u,v)$ replaced by  $a^*(u,v)$. This in turn follows from relations
$\frac{1}{\sqrt{\pi j^*_u}} \e^{-v^2/4j^*_u  } \sim \frac{1}{\sqrt{\pi u(1-\theta)/2}} \e^{-v^2/2 u(1-\theta)  }, \
\psi_u(-d) \sim \Gamma(d)^{-1} u^{d-1}, $ and
\begin{eqnarray} \label{binsum}
\sum_{j\in {\cal D}^+(u,v)} {\rm bin}(u-2j, u; \theta) \ \to \  1/2 \quad \text{as } \quad u \to \infty,
\end{eqnarray}
each of which hold uniformly in $0\le v \le u^{5/9}$.
Let us check \eqref{binsum} for instance. Since
$ c_- u < j < c_+ u,  j\in {\cal D}^+(u,v)$ for some $c_\pm > 0$, see above, so $u^{5/9} = o(j^{3/5})$ and
\eqref{binsum}  follows from
\begin{eqnarray} \label{binsum1}
B'(u)  \  \to \ 1/2 \quad \text{and} \quad B''(u) \to  0,
\end{eqnarray}
where $B'(u) := \sum_{j=0}^u {\rm bin}(u-j, u; \theta)\1 (j \, \text{is even}), \,
B''(u) := \sum_{j=0}^{u} {\rm bin}(u-j, u; \theta)\1 ( |j-u(1-\theta)|\ge u^{3/5}) $.
Here, the first relation in \eqref{binsum1} is obvious by well-known properties of binomial coefficients (??) while the second one
follows from \eqref{hoef} according to which $B''(u) \le  C \e^{-2 u^{1/5}} \to  0.  $  This proves
the first relation in \eqref{Jlim}.

The proof of the second relation in \eqref{Jlim} uses Hoeffding's inequality in \eqref{hoef} in a similar way. We have
$a^-(u,v) \le a_1^-(u,v) + a_2^-(u,v)$, where $a^-_1(u,v) := \psi_u(-d) \sum_{0 \le j \le u:   |j -  u(1-\theta)| > u^{3/5} }   {\rm bin}(u-j, u; \theta) \le C u^{d-1} \e^{-2u^{1/5}} $ implying
$\varrho^{3/2-d} J^-_1(u,v) \le C u^{(10/9)(3/2 -d) + (d-1)} \e^{-2 u^{1/5}} \to 0 \, (u \to \infty)$
uniformly in $|v| \le u^{5/9}$. Finally,
\begin{eqnarray*}
a^-_2(u,v)&:=&\psi_u(-d) \sum_{0 \le j \le u:   |j -  u(1-\theta)| \le u^{3/5}, v \ge j^{3/5} }   {\rm bin}(u-j, u; \theta)
{\rm bin} ( (v+j)/2, j; 1/2) \\
&\le&Cu^{d-1}   \sum_{c_1 u  \le j \le u,  v \ge j^{3/5} } \e^{-2 v^2/j} \ \le \ C u^d \e^{-c_2 u^{1/5}}
\end{eqnarray*}
for some positive constants $c_1, c_2 >0$, implying  $\varrho^{3/2-d} a^-_2(u,v) \le C u^{(10/9)(3/2 -d) + d}
\e^{-c_2 u^{1/5}} \to 0 \, (u \to \infty)$
uniformly in $|v| \le u^{5/9}$ as above. This proves \eqref{Jlim} and Proposition  \ref{propex2}, too.
\hfill $\Box$

\medskip

\noi {\it Proof of Proposition \ref{prop1}.} With the notation $\varrho := \rho(t,s)$ we have that $\{(t,s)\in \R^2, s \ge 0 \} \ni (t,s) \mapsto (\varrho, t/\varrho) \in [0,\infty) \times [-1,1] $ is  a 1-1 mapping. Particularly, if  $\varpi = 1$ then
$(\varrho, \arccos (t/\varrho))$ are the polar coordinates of $(t,s) \in \R^2, s \ge 0$.
We use the inequality:
\begin{equation}\label{tri}
\rho(t_1+t_2, s_1+s_2) \le C_\varpi \sum_{i=1}^2 \rho(t_i,s_i), \quad (t_i,s_i)\in \R^2, \ i=1,2,
\end{equation}
which follows from $\rho (t_1 + t_2, s_1 + s_2)^{1 \wedge \varpi} \le \sum_{i=1}^2 \rho (t_i,s_i)^{1 \wedge \varpi} $
with $C_\varpi := 1 \vee 2^{1/\varpi - 1}.$

\smallskip

\noi (i) W.l.g., let $\delta = 1$. Then $\int_{B_1(0,0)} \rho(t,s)^{-h} \d t \d s
\le 4 \int_0^1 u^{\varpi - h} \d u \int_0^{1/u^\varpi} (1+ v^{2/\varpi})^{-h/2} \d v$, where
the inner integral  $=O(1)$ if $h > \varpi$, $=O(u^{h-\varpi})$ if $h < \varpi$, $= O(|\log u|) $ if
$h = \varpi $, as $u \to 0$.
This proves \eqref{rrho1} and \eqref{rrho2} follows analogously.

\smallskip

\noi (ii) After the change of variables:
$u \to \varrho u$, $v \to \varrho^\varpi v$, $\varrho := \rho(t,s)$, we get
\begin{eqnarray}
(\rho^{-h_1} \star \rho^{-h_2})(t,s)
&=&\varrho^{1+\varpi - h_1 - h_2}
\int_{\R^2} \rho(u,v)^{-h_1}
\rho((t/\varrho) + u, (s/\varrho^\varpi) + v)^{-h_2} \d u \d v,
\nn \\
&=&\varrho^{1+\varpi - h_1 - h_2}(I_1 + I_2 + I_{12}),
\label{rho2}
\end{eqnarray}
where
\begin{eqnarray*}
I_1 &:=& \int_{B_\delta(0,0)} \rho(u,v)^{-h_1}
\rho((t/\varrho) + u, (s/\varrho^\varpi) + v)^{-h_2} \d u \d v, \\
I_2 &:=& \int_{B_\delta(-t/\varrho, - s/\varrho^\varpi)} \cdots \, \d u \d v,\qquad
I_{12} \ := \ \int_{B_\delta^c (0,0) \cap B_\delta^c (-t/\varrho, - s/\varrho^\varpi) } \cdots \, \d u \d v
\end{eqnarray*}
with $\delta > 0$ such that $B_\delta (0,0) \cap B_\delta(-t/\varrho, - s/\varrho^\varpi) = \emptyset $ for any $(t,s)\ne (0,0)$.
The integrals
$I_i \le C, i=1,2 $ by \eqref{rrho1} and  $0< h_i < 1 + \varpi, \, i=1,2$. Next,
by H\"older's inequality with $h := h_1 + h_2$,
\begin{eqnarray} \label{I12}
I_{12} \ \le \ \int_{B^c_\delta(0,0)}
\rho(u,v)^{-h} \d u \d v \ \le \ C,
\end{eqnarray}
in view of \eqref{rrho2} and
$\int_{B^c_\delta (-t/\varrho, -s/\varrho^{\varpi}) } \rho((t/\varrho)+u,(s/\varrho^\varpi)+v)^{-h} \d u \d v =  \int_{B^c_\delta(0,0)}
\rho(u,v)^{-h} \d u \d v. $
This proves \eqref{prop11}.

Next, consider \eqref{prop12}, or
the case $0 < h_2 < 1+ \varpi < h_1$. By changing the variables as in \eqref{rho2}, we get
$(\rho_+^{-h_1} \star \rho^{-h_2})(t,s) \le \varrho^{1+\varpi - h_1 - h_2} (I'_1 + I_2 + I_{12})$,
where $I_2 < C$, $I_{12} < C$ are the same as in \eqref{rho2}, whereas
\begin{equation*}
I'_1 \ := \ \int_{B_\delta (0,0)} (\varrho^{-1} \vee \rho(u,v))^{-h_1} \rho ( (t/\varrho) + u, (s/\varrho^\varpi) + v )^{-h_2} \d u \d v.
\end{equation*}
Note that if given small enough $\delta > 0$, then \eqref{tri} implies $\rho((t/\varrho)+u,(s/\varrho^\varpi)+v)^{1 \wedge \varpi} \ge 1 - \rho(u,v)^{1 \wedge \varpi} \ge 1/2$ for all $(u,v) \in B_\delta (0,0)$, and hence
$I'_1 \le C \varrho^{h_1-1-\varpi} \int_{\R^2} \rho_+ (u,v)^{-h_1} \d u \d v \le C \varrho^{h_1-1-\varpi}$  according to \eqref{rrho2}.
Since $\rho(t,s)^{1+\varpi-h_1-h_2} =  o(\rho(t,s)^{- h_2}) $ as $|t|+|s| \to \infty$, the proof of \eqref{prop12} is complete.

Finally, consider \eqref{prop13}. We follow the proof of \eqref{prop12} and get $(\rho_+^{-h_1} \star \rho_+^{-h_2}) (t,s) \le \varrho^{1+\varpi-h_1-h_2} (I'_1 + I'_2+I_{12})$ with the same $I'_1 < C,$ $I_{12} < C$, whereas
$$
I'_2 \ := \ \int_{B_\delta (-t/\varrho, -s/\varrho^\varpi)} \rho(u,v)^{-h_1} (\varrho^{-1} \vee \rho ( (t/\varrho) + u, (s/\varrho^\varpi) + v ))^{-h_2} \d u \d v.
$$
For small enough $\delta > 0$, we have $\rho(u,v)^{1 \wedge \varpi} \ge 1 - \rho ( (t/\varrho) + u, (s/\varrho^\varpi) + v )^{1 \wedge \varpi} \ge 1/2$ for all $(u,v) \in B_\delta (-t/\varrho,-s/\varrho^\varpi)$, and hence
$I'_2 \le C \varrho^{h_2-1-\varpi} \int_{\R^2} \rho_+ ( (t/\varrho) + u, (s/\varrho^\varpi) + v )^{-h_2} \d u \d v \le C \varrho^{h_2-1-\varpi}$ by~\eqref{rrho2}.
Using $\rho(t,s)^{1+\varpi-h_1-h_2} =  o(\rho(t,s)^{- h_1 \wedge h_2}) $ as $|t|+|s| \to \infty$, we conclude \eqref{prop13}.
Extension of \eqref{prop11}-\eqref{prop13} to
`discrete' convolution $[\rho_+^{-h_1} \star \rho_+^{-h_2}](t,s)$ requires minor changes and
we omit the details. This proves part~(ii).

\smallskip

\noi (iii) It suffices to show  \eqref{prop14} for $(t,s) \ne (0,0),  s \ge 0$, in which case
$\rho_+(t,s)  = \rho(t,s)$. We have
$[a_1 \star a_2](t,s) = \sum_{i,j=0}^1 [a^i_1 \star a^j_2](t,s),$
where $a^0_i(t,s) := \rho_+(t,s)^{-h_i}L_i(t/\rho_+(t,s)), a^1_i(t,s) := a_i(t,s) - a^0_i(t,s) = o(\rho_+(t,s)^{-h_i}), \, i=1,2 $.
Clearly, \eqref{prop14} follows from
\begin{equation}\label{aa0}
\lim_{|t|+|s| \to \infty} \big| \rho(t,s)^{h_1 + h_2 - 1 - \varpi}[a^0_1 \star a^0_2](t,s) - L_{12}(t/\rho(t,s))\big| \ = \  0
\end{equation}
and
\begin{equation}\label{aa1}
 [a^i_1 \star a^j_2](t,s)  \ = \ o(\rho(t,s)^{1+ \varpi - h_1 -  h_2}), \quad (i,j) \ne (0,0), \ i,j =0,1,  \quad |t|+|s| \to \infty.
\end{equation}
The proof of \eqref{aa1} mimics the proof of \eqref{aa0} and is omitted.
To prove  \eqref{aa0}, write $[a^0_1 \star a^0_2](t,s)$ as
the integral:
$[a^0_1 \star a^0_2](t,s)
 = \int_{\R^2} a^0_1([u],[v]) a^0_2([u]+t,[v]+s)  \d u \d v$.
$[a^0_1 \star a^0_2](t,s)
 = \int_{\R^2} a^0_1([u],[v]) a^0_2([u]+t,[v]+s)  \d u \d v$.
After the same change of variables $u \to \varrho u$, $v \to \varrho^\varpi v$, $\varrho := \rho(t,s)$ as in the proof of (ii)
we obtain
$[a^0_1 \star a^0_2](t,s)= \varrho^{1+ \varpi - h_1 - h_2} L_\varrho (t/\varrho),  $
where
\begin{equation}\label{Lrho}
L_\varrho (z) \ := \ \int_{\R^2} g_\varrho (u,v;z) \d u \d v, \quad z \in [-1,1]
\end{equation}
and where
\begin{eqnarray}\label{def:g}
g_\varrho (u,v;z) \ := \ a_{1 \varrho} \big( \tilde u, \tilde v \big) a_{2 \varrho} \big( \tilde u + z, \tilde v + (1-z^2)^{\varpi/2} \big),
\end{eqnarray}
with $\tilde u := [\varrho u]/\varrho$, $\tilde v := [\varrho^\varpi v]/\varrho^\varpi$ and
\begin{eqnarray}
a_{i \varrho} (u, v) \ := \
\big( \varrho^{-1} \vee \rho (u,v) \big)^{-h_i} L_i \big(u/ \big( \varrho^{-1} \vee
\rho (u,v) \big) \big),
\quad i = 1, 2,\label{arrho}
\end{eqnarray}
since $s /\varrho^\varpi = (1 - z^2)^{\varpi/2}$ for $z = t/\varrho  \in [-1,1], s \ge 0$.
Then with $a_{i \infty} (u,v), i=1,2 $ defined by the statement of Prop.~\ref{prop1}~(iii)  we get that
\begin{eqnarray}\label{gconv}
g_\varrho (u,v;z) \ \to \ g_\infty (u,v;z) \ := \ a_{1\infty} (u,v) a_{2\infty} (u+z, v + (1-z^2)^{\varpi/2})
\end{eqnarray}
as $\varrho = \rho (t,s) \to \infty$ ($|t| + |s| \to \infty$)
for any fixed $(u,v; z) \in \R^2 \times [-1,1]$ such that $(u,v) \not \in \{ (0,0), (-z,  -(1-z^2)^{\varpi/2} ) \}$ and $u/\rho (u,v)$, $(u+z)/\rho(u+z, v+(1-z^2)^{\varpi/2})$ being continuity points of $L_1$ and $L_2$ respectively.
Let us prove that
\begin{equation}\label{Llim}
L_\varrho (z) \ \to \   L_{12}(z) \quad \text{as } \varrho \to \infty
\end{equation}
uniformly in $z \in [-1,1]$, which implies \eqref{aa0}, viz., $|L_\varrho (t/\varrho) - L_{12}(t/\varrho)| \le \sup_{z \in [-1,1]} |L_\varrho (z) - L_{12} (z)| = o(1)$ as $\varrho \to \infty$. The uniform convergence in \eqref{Llim} follows if $\lim_{\varrho \to \infty} L_\varrho (z_\varrho) = L_{12}(z)$ holds for any $z \in [-1,1]$ and every sequence $\{ z_\varrho \} \subset [-1,1]$ tending to $z$: $\lim_{\varrho \to  \infty } z_\varrho
= z$. Choose $\delta > 0$ and split the difference
$L_\varrho (z_\varrho ) - L_{12}(z) = I_1  + I_2 + I_{12},$
where
\begin{eqnarray*}
I_1&:=& \int_{B_\delta (0,0)} (g_\varrho(u,v;z_\varrho) - g_\infty(u,v;z)) \d u \d v,\\
I_2 &:=& \int_{B_\delta (-z,-z')} \cdots \, \d u \d v,\qquad
I_{12} \ := \ \int_{B^c_\delta (0,0) \cap B_\delta^c (-z,-z')} \cdots \, \d u \d v
\end{eqnarray*}
with the notation $z' := (1-z^2)^{\varpi/2}$. Note that $\rho (z,z') = 1$ and
$\delta>0$ is chosen small enough so that $B_\delta (0,0) \cap B_\delta (-z,-z') = \emptyset$.
Let us first check that $|I_i|, i=1,2 $ can be made arbitrary small by taking sufficiently small $\delta$.
Towards this end, we need the bound
\begin{equation}\label{a0bdd}
|a_{i \varrho} (\tilde u, \tilde v)|\ \le\  C\rho (u,v)^{-h_i}, \quad (u,v) \in \R^2,  \quad i = 1, 2.
\end{equation}
Indeed, by \eqref{tri}, $\rho(u,v) \le  C_\varpi (\rho(\tilde u, \tilde v) + \rho(u-\tilde u, v - \tilde v))$, where
$|u-\tilde u| \le \varrho^{-1}, |v- \tilde v| \le \varrho^{-\varpi} $ and  hence $\rho(u-\tilde u, v-\tilde v)
\le  \sqrt{2} \varrho^{-1}, $ with $C_\varpi >0$ dependent only on $\varpi >0$.  Therefore,
$\rho(u,v) \le  \sqrt{2} C_\varpi (\rho(\tilde u, \tilde v) + \varrho^{-1})  \le 2 \sqrt{2} C_\varpi (\rho (\tilde u, \tilde v) \vee \varrho^{-1})$
implying $\rho (\tilde u, \tilde v) \vee \varrho^{-1} \ge (2 \sqrt{2} C_\varpi)^{-1} \rho(u,v)$,
or \eqref{a0bdd} in view of the definition of $a_{i \varrho } $ in \eqref{arrho}.
Using \eqref{a0bdd} it follows that
Using \eqref{a0bdd} it follows that
\begin{equation}\label{ggbdd}
|g_\varrho(u,v;z_\varrho) - g_\infty(u,v;z)| \le C\rho (u,v)^{-h_1}
\big(\rho (u+z_\varrho,v+z_\varrho')^{-h_2} + \rho (u+z,v+z')^{-h_2}\big).
\end{equation}
From \eqref{ggbdd} we obtain
$|I_1| \le C\int_{B_\delta (0,0)} \rho (u,v)^{-h_1} \d u \d v \le C \delta^{1 + \varpi - h_1} = o(1)$
and similarly, $|I_2| \le C \delta^{1 + \varpi - h_2} = o(1). $
Hence it suffices to show that
$I_{12} \to 0 \, (z_\varrho \to z)$,
viz., that for each $\delta >0$
\begin{equation}\label{I12z}
\int_{B^c_\delta (0,0) \cap B_\delta^c (-z,-z')} | g_\varrho(u,v;z_\varrho) - g_\infty(u,v;z) | \d u \d v \ \to \ 0 \quad \text{as } \varrho \to \infty.
\end{equation}
From \eqref{tri}, $\rho (u+ z_\varrho, v + z'_\varrho)^{1 \wedge \varpi} \ge \rho(u+z, v+z')^{1 \wedge \varpi} - \delta^{1 \wedge \varpi}/2 \ge (1/2) \rho(u+z, v+z')^{1 \wedge \varpi}$ for all $(u,v) \in  B_\delta^c (-z,-z')$ and $\varrho$ large enough that
$ \rho (z - z_\varrho, z' - z'_\varrho)^{1 \wedge \varpi} \le \delta^{1 \wedge \varpi}/2$ (in view of $z_\varrho \to z$).  Hence and from
\eqref{ggbdd} we obtain that the integrand in \eqref{I12z} is dominated on $B^c_\delta (0,0) \cap B_\delta^c (-z,-z')$
by an integrable function independent of $\varrho$, viz.,
$|g_\varrho(u,v;z_\varrho) - g_\infty(u,v;z)| \le C\rho (u,v)^{-h_1} \rho (u+ z, v + z')^{-h_2}$. Since this integrand
vanishes a.e.\ on  $B^c_\delta (0,0) \cap B_\delta^c (-z,-z')$ as $\varrho \to \infty$, see \eqref{gconv}, relation \eqref{I12z} follows by the dominated
convergence theorem, proving \eqref{Llim}. The continuity of $L_{12}$ \eqref{L12} follows similarly by the dominated
convergence theorem.

It remains to prove the strict positivity of $L_{12}$ in the case where  $L_1(z) \equiv L_2(z) =: L(z)\ge 0.$ Under assumption of piecewise continuity of $L$ and $L \not \equiv 0$ a.e., we can find $0 <|z_0| <1$ and $\delta >0$ such that $L(z)> \delta$ for any $|z- z_0| < \delta$. We also have
$| u/\rho(u,v) - (u+z)/\rho(u+z,v+z') | \le \rho(u,v)^{-1} +
 | 1 - \rho(u+z,v+z')/\rho(u,v) | = O(\rho(u,v)^{-1 \wedge \varpi}) $
uniformly in $z \in [-1,1]$ for $\rho(u,v) \ge 1$.
Indeed, this follows from $|1 - (\rho(u+z,v+z')/\rho(u,v))^{1 \wedge \varpi}| < \rho(u,v)^{-1 \wedge \varpi}$ by \eqref{tri}, when combined with $|1-x| \le \varpi^{-1} (1 \vee x)^{1-\varpi} | 1-x^\varpi |$, $x > 0$ due to the mean value theorem if $0 < \varpi < 1$.
Hence for $\rho(u,v) \ge K \ge 1$ large enough, we have $|u/\rho(u,v) - (u+z)/\rho(u+z,v+z')| < \delta/2$. Next, we can find the interior point $(u_0,v_0)$ of $B_K^c (0,0)$ with $u_0/\rho(u_0,v_0) = z_0$. In view of continuity of $u/\rho(u,v)$, there exists $\varepsilon > 0$ such that $|z_0 - u/\rho(u,v)| < \delta/2$ holds for all $(u,v) \in B_\varepsilon (u_0, v_0) \subset B^c_K (0,0)$. Consequently, $L(u/\rho(u,v))L((u+z)/\rho(u+z,v+z')) \ge \delta^2 >0$ for any $z \in [-1,1]$ and all $(u,v) \in B_\varepsilon (u_0, v_0)$. By \eqref{tri}, we have $\rho(u+z,v+z') \le 2 C_\varpi \rho (u,v)$ for $\rho(u,v) \ge 1$ and hence $L_{12}(z) > \delta^2 (2 C_\varpi)^{-h_2}  \int_{B_\varepsilon (u_0, v_0)} \rho(u,v)^{-h_1-h_2}  \d u \d v  >0,$
proving $L_{12}(z) >0, z \in [-1,1]$ and part (iii).

\smallskip

\noi (iv) Rewrite the l.h.s. of \eqref{Bconv}  as
\begin{equation}\label{Bint}
B_\lambda(\gamma) \ = \ \int_{\widetilde K^2_{[\lambda, \lambda^\gamma]}} b([t_1]-[t_2], [s_1]-[s_2]) \d t_1 \d t_2 \d s_1 \d s_2,
\end{equation}
where $\widetilde K_{[\lambda, \lambda^\gamma]} := \{ (t,s) \in \R^2: ([t],[s]) \in K_{[\lambda, \lambda^\gamma]}\} $.

\smallskip

\noi Case (I): $\gamma = \varpi$. By
changing the variables in \eqref{Bint} as $t_i \to \lambda t_i, s_i \to \lambda^{\varpi} s_i, \, i=1,2,$ we obtain
$\lambda^{-2 {\cal H} (\varpi)}B_{\lambda}(\varpi) = \int_{\R^4} \widetilde b_\lambda(t_1,t_2,s_1,s_2) \d t_1 \d t_2 \d s_1 \d s_2 $, where
\begin{eqnarray}\label{Glambda}
\widetilde b_\lambda(t_1,t_2,s_1,s_2) &:=& b_\lambda ( ([\lambda t_1] - [\lambda t_2])/\lambda, ([\lambda^\varpi s_1]- [\lambda^\varpi s_2])/\lambda^\varpi)\\
&&\times \1 ( ([\lambda t_i], [\lambda^\varpi s_i]) \in (0, \lambda] \times (0, \lambda^\varpi], i=1,2 )\nn
\end{eqnarray}
with $b_\lambda (t,s) := (\lambda^{-1} \vee \rho (t,s) )^{-h} ( L(t/(\lambda^{-1} \vee \rho(t,s)) + o(1) )$ as $\lambda \to \infty$.
Then
\begin{eqnarray*}
\widetilde b_\lambda(t_1,t_2,s_1,s_2) \ \to \ b_\infty (t_1-t_2, s_1-s_2) \1( (t_i,s_i) \in (0,1]^2, i=1,2),  \qquad  \lambda \to \infty
\end{eqnarray*}
point-wise for any $(t_1,t_2,s_1,s_2)\in \R^4, (t_1,s_1) \ne (t_2,s_2) $ fixed.
The dominating bound
\begin{eqnarray*}
&\lambda^{-1} \vee \rho \big( ([\lambda t_1] - [\lambda t_2])/\lambda, ([\lambda^\varpi s_1]- [\lambda^\varpi s_2])/\lambda^\varpi \big) \ \ge \ C \rho (t_1-t_2, s_1-s_2),
\end{eqnarray*}
follows by the same arguments as \eqref{a0bdd}.
These facts and the dominated convergence theorem justify the limit
$\lim_{\lambda \to \infty} \lambda^{-2 {\cal H} (\varpi)}B_{\lambda}(\varpi) = {\cal C} (\varpi)$ since
the integral ${\cal C} (\varpi) \le C \int_{(-1,1]^2} \rho(t,s)^{-h} \d t \d s < \infty$ in \eqref{HC}
converges by Prop.~\ref{prop1}~(i).

\smallskip

\noi Case (II):  $\gamma > \varpi,  h < \varpi$. By changing the variables in \eqref{Bint} as
$t_i \to \lambda t_i,  \ s_i \to \lambda^\gamma s_i, \ i=1,2,$
we obtain $\lambda^{-2{\cal H}(\gamma)}B_{\lambda}(\gamma) = \int_{\R^4} \widetilde b_\lambda(t_1,t_2,s_1,s_2)  \d t_1 \d t_2 \d s_1 \d s_2 $, where
\begin{eqnarray*}
\widetilde b_\lambda(t_1,t_2,s_1,s_2)
&:=&b_\lambda ( ([\lambda t_1] -
[\lambda t_2])/\lambda^{\gamma/\varpi}, ([\lambda^\gamma s_1]-[\lambda^\gamma s_2])/\lambda^\gamma )\\
&&\times \1 ( ([\lambda t_i], [\lambda^\gamma s_i]) \in (0, \lambda] \times (0, \lambda^\gamma], i=1,2 )
\end{eqnarray*}
with $b_\lambda (t,s) := ( \lambda^{-\gamma/\varpi} \vee \rho(t,s) )^{-h} (L(t/(\lambda^{-\gamma/\varpi} \vee \rho(t,s))) + o(1))$ as $\lambda \to \infty$. Hence since $\gamma/\varpi > 1$ it follows that
\begin{eqnarray*}
\widetilde b_\lambda(t_1,t_2,s_1,s_2) \ \to \ b_\infty (0, s_1-s_2) \1( (t_i,s_i) \in (0,1]^2, i=1,2), \qquad  \lambda \to \infty
\end{eqnarray*}
point-wise for any $(t_1,t_2,s_1,s_2)\in \R^4, s_1 \ne s_2$ fixed.  Note $b_\infty (0, s) = L(0)|s|^{-h/\varpi} $ is integrable on $[-1,1]$ due to $h < \varpi$.
The limit $\lim_{\lambda \to \infty} \lambda^{-2 {\cal H} (\varpi)}B_{\lambda}(\varpi) = {\cal C} (\varpi)$  can be justified by the dominated convergence theorem using the bound
\begin{eqnarray*}
\lambda^{-\gamma/\varpi} \vee \rho \big( ([\lambda t_1] - [\lambda t_2])/\lambda^{\gamma/\varpi}, ([\lambda^\gamma s_1]-[\lambda^\gamma s_2])/\lambda^\gamma \big) &\ge& \lambda^{-\gamma/\varpi} \vee \rho ( 0, ([\lambda^\gamma s_1]-[\lambda^\gamma s_2])/\lambda^\gamma )\\
&\ge& C \rho (0, s_1-s_2),
\end{eqnarray*}
which follows by the same arguments as \eqref{a0bdd}.

\smallskip

\noi Case (III):  $\gamma > \varpi,  h > \varpi$.  By changing the variables in \eqref{Bint} as $t_i \to \lambda t_i, i=1,2, s_1- s_2 \to \lambda^{\varpi} s_1, s_2 \to \lambda^\gamma s_2,$
we obtain
$\lambda^{-2{\cal H}(\gamma)}B_\lambda (\gamma) = \int_{\R^4} \widetilde b_\lambda (t_1,t_2,s_1,s_2)  \d t_1 \d t_2 \d s_1 \d s_2, $
where
\begin{eqnarray*}
\widetilde b_\lambda (t_1,t_2,s_1,s_2) &:=& b_\lambda (([\lambda t_1] -
[\lambda t_2])/\lambda, ([\lambda^\varpi s_1 + \lambda^\gamma s_2  ]-[\lambda^\gamma s_2])/\lambda^\varpi )\\
&&\times \1 ( [\lambda t_i] \in (0, \lambda], i=1,2, [\lambda^\varpi s_1 + \lambda^\gamma s_2] \in (0, \lambda^\gamma],
[\lambda^\gamma s_2] \in (0, \lambda^\gamma] )
\end{eqnarray*}
with $b_\lambda (t,s) := (\lambda^{-1} \vee \rho (t,s) )^{-h} ( L(t/(\lambda^{-1} \vee \rho(t,s)) + o(1) )$ as $\lambda \to \infty$. Then
\begin{eqnarray*}
\widetilde b_\lambda(t_1,t_2,s,u) & \to & b_\infty (t_1-t_2, s_1) \1( (t_1,t_2,s_2) \in (0,1]^3), \qquad  \lambda \to \infty
\end{eqnarray*}
for any $t_1 \ne t_2, s_1 \in \R \setminus \{0\}, s_2 \in \R \setminus \{0,1\} $ fixed since $\gamma > \varpi$ implies
$ \1 (0 < [\lambda^\varpi s_1 + \lambda^\gamma s_2] \le \lambda^\gamma]) \to \1 (0< s_2 < 1)$.
The dominating bound
$$
\lambda^{-1} \vee \rho \big( ( [\lambda t_1] - [\lambda t_2] )/ \lambda, ( [\lambda^\varpi s_1 + \lambda^\gamma s_2] - [\lambda^\gamma s_2] )/\lambda^\varpi \big) \ \ge \ C \rho(t_1 - t_2, s_1)
$$
follows by the same arguments as \eqref{a0bdd}, because $|( [\lambda^\varpi s_1 + \lambda^\gamma s_2] - [\lambda^\gamma s_2] )/\lambda^\varpi - s_1 | \le 2 \lambda^{-\varpi}$. Then the dominated convergence in \eqref{Bconv} is proved in view of ${\cal C} (\gamma) \le C \int_{-1}^1 \int_\R \rho(t,s)^{-h} \d t \d s < \infty.$

\smallskip

\noi Cases (IV) and (V) can be treated similarly to Cases (II) and (III) and we omit the details.  Proposition \ref{prop1}
is proved.  \hfill $\Box$

\medskip

\noi {\it Proof of Proposition \ref{appell}.}  (i) We first prove \eqref{Zcov}. According
to \eqref{Awick} 
\begin{eqnarray}\label{Zap}
\hskip-.5cm Z(t,s)&=&\sum_{i=1}^{k-1} \sum_{(D)_i} {\sum_{(u,v)_i}}^{\!\!\!\bf \bullet}
a(t-u_1,s-v_1)^{|D_1|} \cdots a(t-u_i,s-v_i)^{|D_i|} A_{|D_1|}(\vep(u_1,v_1)) \cdots A_{|D_i|}(\vep(u_i,v_i))
\end{eqnarray}
where the sum $\sum_{(D)_i} $ is taken over all partitions of $\{1, 2, \cdots, k \}$ into $i$ nonempty sets
$D_1, \cdots, D_i $ having cardinality $|D_1|\ge 1, \cdots, |D_i| \ge 1, \, |D_1| + \cdots + |D_i| = k$.
Thus,
\eqref{Zap} is a decomposition of $Z(t,s) =  A_k(Y(t,s)) - Y^{\bullet k}(t,s)$ into a sum of stationary
`off-diagonal' polynomial forms of order $i <k$
in i.i.d. r.v. $A_{|D_\ell|}(\vep(u_\ell,v_\ell)), 1\le \ell \le i$ with
$\max (|D_1|, \cdots, |D_i|) \ge 2 $.
From \eqref{Zap} it follows that
\begin{eqnarray}\label{Zap1}
|\E Z(0,0) Z(t,s)|&\le&C \sum_{i=1}^{k-1} \sum_{(d)_i, (d')_i}  \prod_{\ell=1}^i (|a|^{d_\ell} \star |a|^{d'_{\ell}})(t,s)
\end{eqnarray}
where the second sum  is taken over all collections $(d)_i = (d_1, \cdots, d_i),  (d')_i = (d'_1, \cdots, d'_i) $ of integers
$d_\ell \ge 1, d'_\ell \ge 1 $ with $\sum_{\ell=1}^i d_\ell = \sum_{\ell =1}^i d'_\ell = k $ and satisfying
$\max_{1\le \ell \le i} d_\ell \ge 2,  \max_{1\le \ell \le i} d'_\ell \ge 2 $. See \cite{gir2012}, proof of Thm 14.2.1.
Then
$a(t,s)^{d_\ell} \le C \rho(t,s)^{-\beta_\ell}, a(t,s)^{d'_\ell} \le C \rho(t,s)^{-\beta'_\ell} $ where $\beta_\ell := d_\ell q_1, \beta'_\ell
:= d'_\ell q_1 $.  By
Proposition \ref{prop1} (ii)
\begin{eqnarray}\label{Zbdd}
|\E Z(0,0) Z(t,s)|&\le&C \sum_{i=1}^{k-1} \sum_{(d)_i, (d')_i}  \prod_{\ell=1}^i \rho(t,s)^{- w_\ell},
\end{eqnarray}
where
\begin{equation}\label{well}
w_\ell \ := \ \begin{cases}
2q_1-1-\gamma_0 = p_1, &\text{if} \ d_\ell = d'_\ell = 1, \\
q_1, &\text{if} \ d_\ell \ge 2, \, d'_\ell =1 \ \text{or} \  d_\ell = 1, \, d'_\ell \ge 2, \\
2q_1, &\text{if} \ d_\ell \ge 2, \, d'_\ell \ge 2.
\end{cases}
\end{equation}
Since $2q_1 - 1 - \gamma_0 < q_1 < 2q_1 $, we have that, for  $\max_{1\le \ell \le i} d_\ell \ge 2,  \max_{1\le \ell \le i} d'_\ell \ge 2 $,
the exponents $w_\ell$ in \eqref{well} satisfy $\sum_{\ell=1}^i w_\ell \ge 2q_1 $, implying \eqref{Zcov}.
Since RFs $\{Y^{\bullet k}(t,s)\}$ and $\{Z(t,s)\}$ are uncorrelated:
${\rm Cov}(Y^{\bullet k}(t,s), Z(u,v)) $  $ =0 $ for any $(t,s), (u,v) \in \Z^2$,
relation \eqref{rhoas2} follows from \eqref{Zcov} and
\begin{equation}\label{Ykcov}
{\rm Cov}(Y^{\bullet k}(t,s),Y^{\bullet k}(0,0)) \ = \  r_Y(t,s)^k(1 + o(1)), \qquad
|t|+ |s| \to  \infty.
\end{equation}
To show \eqref{Ykcov}, note that the difference $|r_Y(t,s)^k -
{\rm Cov}(Y^{\bullet k}(t,s),Y^{\bullet k}(0,0))| =  \big|(a \star a)^k(t,s) - \sum_{(u,v)_k}^{\bf \bullet}  a(t+u_1,s+v_1) a(u_1,v_1)
 \cdots a(t+u_k,s+v_k)  a(u_k,v_k)\big|  $ satisfies the same bound as in \eqref{Zbdd} 
 and
therefore this difference is   $O(\rho(t,s)^{-2q_1}) =  o(r_Y(t,s)^k) $ according to \eqref{Zcov}.
This proves \eqref{Ykcov} and part (i). Part (ii) follows similarly using \eqref{Zcov} and $|{\rm Cov}(Y^{\bullet k}(t,s),Y^{\bullet k}(0,0))| \le (|a| \star |a|(t,s))^k
\le C \rho_+(t,s)^{-kp_1}$.
Proposition \ref{appell} is proved. \hfill $\Box$

\medskip

\noi {\it Proof of Corollary \ref{appellvar}.} 
Relation \eqref{varZ} follows from  \eqref{Zcov} and Proposition \ref{prop1}~(i) since the l.h.s. of
\eqref{varZ} does not exceed
$ \sum_{(t_1,s_1), (t_2,s_2) \in K_{[\lambda,\lambda^{\gamma}]}}
|r_Z (t_1-t_2,s_1-s_2)| \le \lambda^{1+ \gamma} \sum_{(t,s)\in \Z^2} |r_Z (t,s)| \le
C \lambda^{1+ \gamma} \sum_{(t,s)\in \Z^2} \rho_+(t,s)^{-2q_1} $ and  the last sum  converges by Proposition \ref{prop1} (i)
due to $2q_1 > 1 + \gamma_0$. \\
Relations 
\eqref{varX} follow from \eqref{varZ},
the orthogonality of $\{  Y^{\bullet k}(t,s)\} $ and  $\{ Z(t,s)\} $ and
\begin{eqnarray}\label{varYk}
{\rm Var}\big(\sum\nolimits_{(t,s) \in K_{[\lambda,\lambda^{\gamma}]}}  Y^{\bullet k}(t,s) \big)
&\sim&c(\gamma)\lambda^{2H(\gamma)}.
\end{eqnarray}
In turn, \eqref{varYk} follows from
\begin{eqnarray}\label{varYk1}
{\cal V}_{\lambda,\gamma} \ := \
\sum_{(t_1,s_1), (t_2,s_2) \in K_{[\lambda,\lambda^{\gamma}]}}
r^k_Y (t_1-t_2,s_1-s_2)&\sim&c(\gamma)^2 \lambda^{2H(\gamma)}
\end{eqnarray}
and the fact that the difference $\big|{\rm Var}\big(\sum\nolimits_{(t,s) \in K_{[\lambda,\lambda^{\gamma}]}}  Y^{\bullet k}(t,s) \big)
- \sum_{(t_1,s_1), (t_2,s_2) \in K_{[\lambda,\lambda^{\gamma}]}}
r^k_Y (t_1-t_2,s_1-s_2)\big| $ can be estimated as in  \eqref{Zap1}-\eqref{Zbdd} and therefore
this difference is $O(\lambda^{1+\gamma}) = o(\lambda^{2H(\gamma)} $ as shown in  \eqref{varZ}.


\vskip.2cm
































\bigskip




\footnotesize


\begin{thebibliography}{99}





\bibitem{AnhLRM2012} Anh, V.V., Leonenko, N.N. and Ruiz-Medina, M.D. (2013) Macroscaling limit theorems for filtered spatiotemporal random fields.
Stochastic Anal. Appl. 31, 460--508.


\bibitem{avr1987} Avram, F. and Taqqu, M.S. (1987) Noncentral limit theorems and Appell polynomials. Ann. Probab. 15, 767--775.

\bibitem{bai2014} Bai, S. and Taqqu, M.S. (2014) Generalized Hermite processes, discrete chaos and
limit theorems. Stochastic Process. Appl. 124, 1710--1739.


\bibitem{bier2007} Bierm\'e, H., Meerschaert, M.M. and Scheffler, H.P. (2007) Operator scaling
stable random fields. Stoch. Process. Appl. 117, 312-332.



\bibitem{bolt1982} Bolthausen, E. (1982) On the central  limit theorem for stationary mixing random
fields. Ann. Probab. 10, 1047--1050.

\bibitem{breu1983} Breuer, P. and Major, P. (1983) Central limit theorems for
non-linear functionals of Gaussian fields.  J. Multiv. Anal. 13, 425-441.



\bibitem{bois2005}  Boissy, Y., Bhattacharyya, B.B., Li, X. and Richardson, G.D. (2005) Parameter estimates
for fractional autoregressive spatial processes. Ann. Statist. 33, 2533--2567.













\bibitem {dobr1979} Dobrushin, R.L. (1979) Gaussian and their subordinated self-similar random
generalized fields. Ann. Probab. 7, 1--28.


\bibitem {dobmaj1979} Dobrushin, R.L. and Major, P. (1979)
Non-central limit theorems for non-linear functionals of
Gaussian fields. Z. Wahrsch. verw. Geb.
50, 27--52.








\bibitem{fell1966} Feller, W. (1966) {\em An Introduction to Probability Theory and Its Applications}, vol. 2. Wiley,
New York.







\bibitem{gaig2003} Gaigalas, R. and Kaj, I. (2003)
Convergence of scaled renewal processes and a packet arrival model.
Bernoulli 9, 671--703.





\bibitem{guo2009} Guo, H., Lim, C. and Meerschaert, M. (2009) Local Whittle estimator for anisotropic random
fields. J. Multiv. Anal. 100, 993--1028.

\bibitem{gir1985} Giraitis, L. and  Surgailis, D.   (1985) CLT and other limit theorems for functionals
of Gaussian processes. Z. Wahrsch. verw. Geb.
70, 191--212.

\bibitem{gir2012} Giraitis, L., Koul, H.L. and  Surgailis, D.   (2012) {\em Large Sample Inference
for Long Memory Processes.} Imperial College Press, London.


\bibitem{han1972} Hankey, A. and Stanley, H.E.  (1972)
Systematic application of generalized homogeneous functions to static
scaling, dynamic scaling, and universality. Phys. Review B 6,
3515--3542.


\bibitem{hoh1997} Ho, H.-C. and Hsing, T. (1997) Limit theorems for functionals of moving averages.
Ann. Probab. 25, 1636--1669.


\bibitem{hoef1963} Hoeffding, W. (1963) Probability inequalities for sums of bounded random variables.
J. Amer. Statist. Assoc. 58, 13--30.




\bibitem{kou2016} Koul, H.L., Mimoto, N. and Surgailis, D. (2016)
Goodness-of-fit tests for marginal distribution of linear
random fields with long memory. Metrika 79, 165--193.



\bibitem{lav2007} Lavancier, F. (2007) Invariance principles for non-isotropic long memory
random fields. Statist. Inference Stoch. Process. 10, 255--282.



\bibitem{lls2014} Lavancier, F., Leipus, R. and Surgailis, D. (2014) Anisotropic long-range dependence and
aggregation of space-time models. Preprint.



\bibitem{leo1999} Leonenko, N.N. (1999) {\em
Random Fields with Singular Spectrum}. Kluwer, Dordrecht.





\bibitem{miko2002} Mikosch, T., Resnick, S., Rootz\'en, H. and Stegeman, A. (2002)
Is network traffic approximated by stable L\'evy motion or fractional Brownian motion?
Ann. Appl. Probab. 12, 23--68.

\bibitem{nev1964} Neveu, J. (1964) {\em Bases math\'ematiques du calcul de probabilit\'es}, Mason, Paris.

\bibitem{nua2005} Nualart, D. and Peccati, G. (2005) Central limit theorems for sequences of
multiple stochastic integrals. Ann. Probab. 33, 177--193.

\bibitem{orey1958} Orey, S. (1958) A central limit theorem for $m$-dependent random variables. Duke Math. J. 25, 543--546.


\bibitem{pils2014} Pilipauskait\.e, V.  and Surgailis, D. (2014) Joint temporal and contemporaneous aggregation
of random-coefficient  AR(1) processes. Stochastic Process. Appl. 124, 1011--1035.


\bibitem{pils2015} Pilipauskait\.e, V.  and Surgailis, D. (2015) Joint aggregation of random-coefficient AR(1) processes
with common innovations. Statist. Probab. Letters 101, 73--82.




\bibitem{pils2016} Pilipauskait\.e, V.  and Surgailis, D. (2016) Anisotropic scaling of random grain model  with application
to network traffic. J. Appl. Probab. (in press). Available at http://arxiv.org/abs/1510.07423.





\bibitem{pra1960} Pratt, J.W. (1960) On interchanging limits and integrals. Ann. Math. Statist.
31, 74--77.


\bibitem{ps2014} Puplinskait\.e, D.  and Surgailis, D. (2016)
Aggregation of autoregressive random
fields and anisotropic long-range dependence. Bernoulli (in press), DOI 10.3150/15-BEJ733.
Available at
http://arxiv.org/abs/1303.2209v3.


\bibitem{ps2015} Puplinskait\.e, D.  and Surgailis, D. (2015)
Scaling  transition for long-range dependent Gaussian random fields.
Stoch. Process. Appl. 125, 2256--2271.


\bibitem{sko1956} Skorokhod, A.V. (1956) Limit theorems for stochastic processes. Th. Probab. Appl. 1, 261--290.


\bibitem {sur1982} Surgailis, D. (1982) Zones of attraction of self-similar multiple
integrals. Lithuanian Math. J.  22, 185--201.


\bibitem{taq1979}
Taqqu, M.S. (1979) Convergence of integrated processes of arbitrary Hermite rank.
Z. Wahrsch. verw. Geb. 50, 53--83.
%

Fractional Brownian motion and long-range dependence.



\end{thebibliography}
\end{document}